\documentclass[12pt]{article}

\usepackage{amsmath,amssymb,amsthm,mathtools,bm}
\usepackage{booktabs,multirow,enumitem,float}
\usepackage{ascmac}
\usepackage{xcolor,mathrsfs}
\usepackage{geometry}
\usepackage[expansion=false]{microtype}
\usepackage{etoolbox}
\usepackage{hyperref}
\makeatletter
\let\manuscript@cite\@cite
\renewcommand{\@cite}[2]{{\upshape\manuscript@cite{#1}{#2}}}
\makeatother

\hypersetup{
  pdftitle={Rarita-Schwinger fields and trans-Sasakian eta-Einstein geometry in three dimensions},
  pdfauthor={Satsuki Matsuno},
  pdfkeywords={Rarita-Schwinger field, spin geometry, trans-Sasakian manifold, eta-Einstein metric},
  colorlinks=true,
  linkcolor=blue,
  urlcolor=blue
}

\allowdisplaybreaks
\setlist[enumerate]{itemsep=2pt,parsep=0pt}
\makeatletter
\renewcommand{\paragraph}{\@startsection{paragraph}{4}{\z@}%
  {0.5ex \@plus .5ex \@minus .2ex}%
  {-.5em}%
  {\normalfont\normalsize\bfseries}}
\makeatother

\newcommand{\dd}{\mathop{}\!\mathrm{d}}
\newcommand{\D}{\mathbb D}
\newcommand{\C}{\mathbb C}
\newcommand{\Q}{\mathbb Q}
\newcommand{\R}{\mathbb R}
\newcommand{\Z}{\mathbb Z}

\newcommand{\Ric}{\operatorname{Ric}}
\newcommand{\rank}{\operatorname{rank}}
\newcommand{\Scal}{\operatorname{Scal}}
\newcommand{\Id}{\operatorname{Id}}

\newcommand{\cR}{\Pi_{\Ric}}

\newcommand{\diver}{\operatorname{div}}

\newcommand{\loc}{\mathrm{loc}}

\newcommand{\Lie}{\mathcal L}

\newcommand{\ii}{\mathrm i}

\DeclareMathOperator{\Log}{Log}
\renewcommand{\Re}{\operatorname{Re}}
\renewcommand{\Im}{\operatorname{Im}}

\theoremstyle{plain}
\newtheorem{theorem}{Theorem}[section]
\newtheorem{proposition}[theorem]{Proposition}
\newtheorem{lemma}[theorem]{Lemma}
\newtheorem{corollary}[theorem]{Corollary}
\newtheorem*{corollaryrestate}{Corollary~\ref{coro:S3 regular dimRS}}
\newtheorem*{theoremrestate}{Theorem~\ref{thm:removable-main}}
\newtheorem*{theoremfixedrestate}{Theorem~\ref{thm:fixed}}

\theoremstyle{definition}
\newtheorem{definition}[theorem]{Definition}

\theoremstyle{remark}

\numberwithin{equation}{section}

\title{Rarita-Schwinger fields and trans-Sasakian $\eta$-Einstein geometry in three dimensions}
\author{
Satsuki Matsuno\thanks{Department of Physics, Graduate School of Science, Osaka Metropolitan University.
E-mail: Satsuki\_Matsuno@omu.ac.jp.
ORCID: 0000-0003-0841-8033.}
}
\date{}

\begin{document}
\setlength{\abovedisplayskip}{3.5pt}
\setlength{\belowdisplayskip}{3.5pt}
\maketitle

\begin{abstract}
We study Rarita-Schwinger fields on three-dimensional Riemannian spin manifolds.
Viewing a Rarita-Schwinger field as a real linear map from the tangent space to the spinor space, we show that every Rarita-Schwinger field that is not identically zero on a connected three-dimensional manifold has either (R2) rank two almost everywhere or (R3) rank three almost everywhere.
We further show that the existence of a Rarita-Schwinger field of type (R2) is closely related to a trans-Sasakian $\eta$-Einstein structure.
We then classify Rarita-Schwinger fields of type (R2) on simply connected and geodesically complete three-dimensional Riemannian spin manifolds.
As a corollary, we show that, for a given positive integer, there exists a metric on $S^3$ such that the real dimension of the space of all Rarita-Schwinger fields is greater than that number.
\end{abstract}



\tableofcontents

\section{Introduction}
\label{sec:introduction}

The Rarita-Schwinger equation was introduced by Rarita and Schwinger as a relativistic wave equation describing particles of spin $3/2$~\cite{RaritaSchwinger1941}.
Rarita-Schwinger (RS) fields appear in the study of higher-spin fields and gravitinos and play an important role in mathematical physics, including supergravity.
In Riemannian spin geometry, the RS operator is defined as the spin $3/2$ component of the twisted Dirac operator and can be viewed as a generalization of the usual Dirac operator~\cite{HommaSemmelmann2019}.

RS fields on Riemannian spin manifolds have been studied in relation to the geometry of the underlying manifolds.
Homma and Semmelmann studied the kernel of the RS operator on compact Riemannian spin manifolds and established results on the existence of RS fields and their relation to harmonic forms on quaternionic K\"ahler manifolds, symmetric spaces, and manifolds with special holonomy~\cite{HommaSemmelmann2019}.
B\"ar and Mazzeo also constructed compact spin manifolds in each dimension $n\geq4$ for which the dimension of the space of RS fields can be arbitrarily large~\cite{BarMazzeo2021}.

For manifolds with special geometric structures, Ohno and Tomihisa identified the space of RS fields on six-dimensional compact strict nearly K\"ahler manifolds with the space of harmonic $3$-forms~\cite{OhnoTomihisa2023}, and Ohno described RS fields on compact nearly parallel $G_2$ manifolds in terms of eigenspaces of differential forms~\cite{Ohno2023G2}.
Nguyen also introduced equations of Seiberg--Witten type using the RS operator, which has thus been studied from the viewpoint of gauge theory as well~\cite{Nguyen2023}.

In three-dimension, Ohno classified left-invariant RS fields on connected and simply connected three-dimensional Lie groups equipped with left-invariant metrics~\cite{Ohno2026LieGroups}.
In particular, in the three-dimensional nonabelian case, only $SU(2)$ equipped with a particular Berger metric admits nonzero left-invariant RS fields, and its space of left-invariant RS fields has complex dimension two.
This metric is not Einstein and gives a concrete example of a non-Einstein metric on a compact homogeneous space that admits RS fields.

In view of previous work in low dimensions, it is natural to study the relation between the existence of RS fields and the underlying geometry in three dimensions without assuming any symmetry of the metric or the fields.
A major problem in this area is to determine all RS fields on an arbitrary three-dimensional spin manifold.
Although a complete solution is difficult, this paper gives a partial solution by considering the problem from the viewpoint of the rank of RS fields.

We now describe the rank of an RS field and the scope of this paper.
Since an RS field $\Psi$ is a spinor-valued $1$-form, it defines a real linear map $\Psi_p:T_pM\to\Sigma_pM$ to the fiber of the spinor bundle $\Sigma M$ at each point $p\in M$.
The real rank of this linear map $\Psi_p$ is called the rank of $\Psi$ at $p$.
Since $\dim M=3$, the rank is at most three.
For an RS field that is not identically zero on a connected manifold, Proposition~\ref{prop:rank-propagation} shows that only the following two cases can occur concerning the sizes of the sets on which each rank occurs, namely, (R2) rank two almost everywhere or (R3) rank three almost everywhere (Table~\ref{tab:rank-types}).

\begin{table}[htbp]
\centering
\begin{tabular}{|c|c|c|c|c|}
\hline
Type & Zero & Rank 1 & Rank 2 & Rank 3 \\
\hline
(R2) & \multirow{2}{*}{Measure zero} & \multirow{2}{*}{$\varnothing$} & a.e. & $\varnothing$ \\
\cline{1-1}\cline{4-5}
(R3) & & & Measure zero & a.e. \\
\hline
\end{tabular}
\caption{Each column specifies a condition on $\Psi_p$, and each entry indicates the size of the set on which that condition holds.}
\label{tab:rank-types}
\end{table}

Since $\ker D_{TM}$ on a compact manifold is finite-dimensional, the space of all RS fields is also finite-dimensional, regardless of whether they are of type (R2) or type (R3).
However, in the noncompact case, as will be shown in Sections~\ref{sec:Rarita-Schwinger fields on trans-Sasakian eta-Einstein manifolds}, \ref{sec:explicit-model-fields}, and~\ref{sec:conclusion}, the space of all RS fields of type (R2) on a manifold of nonconstant curvature is infinite-dimensional whenever such fields exist, whereas the space of all RS fields of type (R3) has dimension two whenever such fields exist.
Thus RS fields of type (R2) has a richer structure in terms of its degrees of freedom.

In this paper, we classify RS fields of type (R2) on simply connected complete three-dimensional spin manifolds.
In the case (R2), the kernels $\ker_{\R}\Psi_p$ define a smooth one-dimensional tangent distribution on $M\setminus\{\Psi=0\}$.
By studying the geometry induced by this one-dimensional distribution, we can achieve our aim.
The main results of this paper are as follows.

\begin{theoremrestate}
Let $(M,g)$ be a connected and simply connected three-dimensional Riemannian spin manifold, and let $\Psi$ be an RS field that is not identically zero and is of type (R2), that is, it has rank two almost everywhere on $M$.
Then $M$ admits a trans-Sasakian $\eta$-Einstein structure $(\phi,\xi,\eta,g)$ such that $\ker_\R\Psi_p=\R\xi_p$ at every point $p$ where $\Psi(p)\ne0$.
The vector field $\xi$ satisfying this condition is unique up to sign.
\end{theoremrestate}

\begin{theoremfixedrestate}
Let $(M,\phi,\xi,\eta,g)$ be a connected three-dimensional trans-Sasakian $\eta$-Einstein manifold which does not have constant curvature.  
Every RS field $\Psi$ that is not identically zero is of type (R2) and satisfies $\ker_{\R}\Psi_p=\R\xi_p$ for $\Psi_p\ne0$.
\end{theoremfixedrestate}

The two theorems above imply that RS fields of types (R2) and (R3) cannot coexist on a three-dimensional spin manifold that does not have constant curvature.
Theorem \ref{thm:removable-main} show that, to classify RS fields of type (R2) in dimension three, it suffices to classify trans-Sasakian $\eta$-Einstein manifolds and determine the RS fields on them.
In Section~\ref{sec:explicit-model-fields}, we discuss the following.
Up to homothety, the connected, simply connected and geodesically complete three-dimensional Riemannian spin manifolds that admit an RS field of type (R2) are the following.
\begin{enumerate}
  \item Proper trans-Sasakian $\eta$-Einstein manifolds that are not Einstein.
  \item Sasakian manifolds diffeomorphic to $\R^3$.
  \item $\beta$-Kenmotsu $\eta$-Einstein manifolds diffeomorphic to $\R^3$.
  \item Trans-Sasakian structures on $\mathbb{E}^3$.
  \item Two types of trans-Sasakian structures on $\mathbb{H}^3$.
  \item Sasakian manifolds diffeomorphic to $S^3$ that satisfy some condition.
\end{enumerate}
The RS fields of type (R2) on these manifolds are given in Section~\ref{sec:explicit-model-fields}.
As a further corollary of these results, we obtain metrics on $S^3$ for which the dimension of the space of all RS fields can be arbitrarily large.
B\"ar and Mazzeo showed that, in dimension \(n\ge4\), for any given positive number, there exists a compact manifold on which the space of all RS fields has dimension greater than that number~\cite[Corollary~4.2]{BarMazzeo2021}.
The following corollary is a three-dimensional version of this result of B\"ar and Mazzeo.
Both constructions fix the dimension, but in our construction the underlying manifold is also fixed as $S^3$.

\begin{corollaryrestate}
Let $(S^3,\phi,\xi,\eta,g)$ be a Sasakian structure obtained by a second type deformation of the regular weighted Sasakian structure on $S^3$ with weights $a_1=a_2=\frac{1}{2(d+3)},\ d\in\mathbb Z_{\ge0}$.
Then the complex dimension of the vector space of all RS fields is $2(d+1)$.
\end{corollaryrestate}

\section{Preliminaries}
\label{sec:preliminaries}

\subsection{Rarita-Schwinger fields and real rank}
\label{subsec:rs-fields-real-rank}

Let $(M^3,g)$ be a Riemannian spin manifold, and let $\Sigma M$ be its spinor bundle.
Write $\gamma:TM\to\operatorname{End}(\Sigma M)$ for Clifford multiplication and assume the convention $ \gamma(X)\gamma(Y)+\gamma(Y)\gamma(X)=-2g(X,Y)\operatorname{Id}_{\Sigma M}$ and $\gamma_1\gamma_2\gamma_3=\Id$.
Throughout this paper, all orthonormal frames are assumed to be positively oriented.
For a local orthonormal frame $\{e_1,e_2,e_3\}$, set $\gamma_i:=\gamma(e_i)$.
Define the Clifford contraction by
\begin{equation}\label{eq:clifford-contraction}
 \Pi:T^*M\otimes\Sigma M\longrightarrow\Sigma M,\qquad
 \Pi(\Psi):=\sum_{i=1}^3\gamma_i\Psi(e_i),
\end{equation}
and call its kernel $\Sigma_{3/2}M:=\ker\Pi\subset T^*M\otimes\Sigma M$ the spin $3/2$ bundle.

Denote the spin connection by $\nabla^\Sigma$, the Levi-Civita connection by $\nabla^g$, and the induced connection on $T^*M\otimes\Sigma M$ by $\nabla$.
Thus,
\begin{equation}\label{eq:vector-spinor-connection}
 (\nabla_X\Psi)(Y)
 =\nabla^\Sigma_X(\Psi(Y))-\Psi(\nabla^g_XY).
\end{equation}
Define the usual Dirac operator $D_\Sigma$ and the twisted Dirac operator $D_{TM}$ by
\begin{equation}\label{eq:dirac-operators}
 D_\Sigma\psi:=\sum_{i=1}^3\gamma_i\nabla^\Sigma_{e_i}\psi,\qquad
 (D_{TM}\Psi)(Y):=\sum_{i=1}^3\gamma_i(\nabla_{e_i}\Psi)(Y).
\end{equation}

\begin{definition}\label{def:rs-field}
A spinor-valued $1$-form $\Psi\in\Gamma(\Sigma_{3/2}M)$ is called a Rarita-Schwinger field (RS field)~\cite{HommaSemmelmann2019,OhnoTomihisa2023,Ohno2026LieGroups} if it satisfies $D_{TM}\Psi=0$.
\end{definition}

At each point $p\in M$, regard $\Psi_p:T_pM\to\Sigma_pM$ as a real linear map and define its real rank by $\operatorname{rank}_\R\Psi_p :=\dim_\R\operatorname{im}\bigl(\Psi_p:T_pM\to\Sigma_pM\bigr)$.

\begin{lemma}\label{lem:rs-real-rank}
For $\Psi\in\Gamma(\Sigma_{3/2}M)$ on a three-dimensional Riemannian spin manifold, $\operatorname{rank}_\R\Psi_p\in\{0,2,3\}$ at every point $p\in M$.
In particular, the real rank is $2$ or $3$ at every point where $\Psi_p\ne0$.
\end{lemma}
\begin{proof}
Since $\dim_\R T_pM=3$, the real rank is at most $3$.
Suppose that the real rank is $1$.
Then $\ker_\R\Psi_p$ has real dimension $2$, so we can choose an orthonormal frame $\{e_1,e_2,e_3\}$ such that $e_2,e_3$ form an orthonormal basis of $\ker_\R\Psi_p$.
Since $\Pi(\Psi)=0$, we have $0=\Pi(\Psi)_p=\sum_{i=1}^3\gamma_i\Psi_p(e_i)=\gamma_1\Psi_p(e_1)$, and hence $\Psi_p(e_1)=0$.
Thus we have $\Psi_p=0$.
This contradicts the assumption that the real rank is $1$, so rank $1$ is impossible.
\end{proof}

Define the exterior covariant derivative of the spin connection by $(\dd^\nabla\Psi)(X,Y)=(\nabla_X\Psi)(Y)-(\nabla_Y\Psi)(X)=\nabla^\Sigma_X(\Psi(Y))-\nabla^\Sigma_Y(\Psi(X))-\Psi(\nabla^g_XY)+\Psi(\nabla^g_YX)$.

\begin{lemma}\label{lem:closed}
On a three-dimensional Riemannian spin manifold $(M,g)$, for $\Psi\in\Gamma(\Sigma_{3/2}M)$, $D_{TM}\Psi=0$ is equivalent to $\dd^\nabla\Psi=0$.
\end{lemma}
\begin{proof}
At an arbitrary point $p$, choose a local orthonormal frame such that $\nabla^g_{e_i}e_j(p)=0$.
At $p$, we have $\sum_{i=1}^3\gamma_i(\dd^\nabla\Psi)(e_i,e_j)
=\sum_{i=1}^3\gamma_i\bigl((\nabla_{e_i}\Psi)(e_j)-(\nabla_{e_j}\Psi)(e_i)\bigr)
=(D_{TM}\Psi)(e_j)-\sum_{i=1}^3\nabla^\Sigma_{e_j}(\gamma_i\Psi(e_i))
=(D_{TM}\Psi)(e_j)$.
Since $p$ is arbitrary, this identity holds at every point.

Thus, it remains to show that the map $\Omega^2(M,\Sigma M)\ni F\mapsto\sum_{i=1}^3\gamma_iF(e_i,\,\cdot\,)\in \Omega^1(M,\Sigma M)$ is injective.
Set $F_1=F(e_2,e_3)$, $F_2=F(e_3,e_1)$, and $F_3=F(e_1,e_2)$.
The condition for $F$ to lie in the kernel is $-\gamma_2F_3+\gamma_3F_2=0,\  \gamma_1F_3-\gamma_3F_1=0,\ -\gamma_1F_2+\gamma_2F_1=0$.
Solving these equations together gives $F=0$.
\end{proof}

An integrability condition of an RS field is as follows.

\begin{lemma}\label{lem:ricci}
On a three-dimensional Riemannian spin manifold, an RS field $\Psi\in\Gamma(\Sigma_{3/2}M)$ satisfies
\begin{equation}\label{eq:ricci-constraint}
 \cR(\Psi):=\sum_{i=1}^3\gamma(\Ric(e_i))\Psi(e_i)=0.
\end{equation}
\end{lemma}
\begin{proof}
Choose an orthonormal frame $\{e_i\}$ such that $\nabla^g_{e_i}e_j=0$ at $p\in M$.
At $p\in M$, we have $0=\Pi(D_{TM}\Psi)=\sum_{i,j}\gamma_i\gamma_j(\nabla_{e_j}\Psi)(e_i)=-\sum_j\gamma_j\left(\sum_i\gamma_i(\nabla_{e_j}\Psi)(e_i)\right)-2\sum_{i,j}\delta_{ij}(\nabla_{e_j}\Psi)(e_i)=-\sum_j\gamma_j\left(\sum_i\gamma_i(\nabla_{e_j}\Psi)(e_i)\right)-2\sum_i(\nabla_{e_i}\Psi)(e_i)=-\sum_j\gamma_j\left(\sum_i\nabla^\Sigma_{e_j}\bigl(\gamma_i\Psi(e_i)\bigr)\right)-2\operatorname{div}\Psi=-\sum_j\gamma_j\nabla^\Sigma_{e_j}(\Pi\Psi)-2\operatorname{div}\Psi=-2\operatorname{div}\Psi$,
thus $\operatorname{div}\Psi=0$ holds on $M$.

We use the well-known curvature contraction identity $\gamma(\Ric(X))=2\sum_{i=1}^3\gamma_i\mathrm{Riem}(e_i,X)$,
then we have
$\operatorname{div}(D_{TM}\Psi)-D_\Sigma(\operatorname{div}\Psi)
=\sum_{i,j}\gamma_j\bigl((\nabla_{e_i}\nabla_{e_j}-\nabla_{e_j}\nabla_{e_i})\Psi\bigr)(e_i)
=\sum_{i,j}\gamma_j\mathrm{Riem}(e_i,e_j)\Psi(e_i)-\sum_{i,j}\gamma_j\Psi(R(e_i,e_j)e_i)
=\sum_i\left(\sum_j\gamma_j\mathrm{Riem}(e_i,e_j)\right)\Psi(e_i)-\sum_j\gamma_j\Psi\left(\sum_iR(e_i,e_j)e_i\right)
=-\frac12\sum_i\gamma\bigl(\operatorname{Ric}^{\sharp}(e_i)\bigr)\Psi(e_i)-\sum_j\gamma_j\Psi\bigl(-\operatorname{Ric}^{\sharp}(e_j)\bigr)
=-\frac12\Pi_{\mathrm{Ric}}\Psi+\sum_j\gamma_j\Psi\bigl(\operatorname{Ric}^{\sharp}(e_j)\bigr)
=-\frac12\Pi_{\mathrm{Ric}}\Psi  +\sum_j\gamma_j\Psi\left(\sum_i\operatorname{Ric}(e_j,e_i)e_i\right)
=-\frac12\Pi_{\mathrm{Ric}}\Psi+\Pi_{\mathrm{Ric}}\Psi=\frac12\Pi_{\mathrm{Ric}}\Psi.$
Therefore we have $\Pi_{\mathrm{Ric}}\Psi=0$.
\end{proof}

\subsection{Trans-Sasakian \texorpdfstring{$\eta$}{eta}-Einstein structures}
\label{subsec:trans-sasakian-eta-einstein}

An almost contact metric structure (ACM structure) on a three-dimensional manifold $M$ is a tuple $(\phi,\xi,\eta,g)$ consisting of an endomorphism field $\phi$, a vector field $\xi$, a $1$-form $\eta$, and a Riemannian metric $g$, satisfying
\begin{equation}\label{eq:acm-structure}
 \begin{aligned}
 \phi^2&=-\Id+\eta\otimes\xi,\quad 
 \eta(\xi)=1,\quad 
 g(\phi X,\phi Y)=g(X,Y)-\eta(X)\eta(Y).
 \end{aligned}
\end{equation}
These conditions imply $\phi\xi=0,\ \eta\circ\phi=0,\  \eta=\xi^\flat$ and $g(\xi,\xi)=1$.

An ACM structure is called \textcolor{blue}{\textbf{trans-Sasakian}} if there exist smooth real functions $\alpha,\beta$ such that
\begin{align}
(\nabla^g_X\phi)Y
 ={}&\alpha\{g(X,Y)\xi-\eta(Y)X\}
     +\beta\{g(\phi X,Y)\xi-\eta(Y)\phi X\}.
 \label{eq:trans-sasakian}
\end{align}
In particular, in three dimensions, this condition is known to be equivalent to
\begin{equation}\label{eq:xi-connection}
 \nabla^g_X\xi=-\alpha\phi X+\beta\{X-\eta(X)\xi\}.
\end{equation}
By~\eqref{eq:trans-sasakian} and~\eqref{eq:xi-connection}, we have $ \nabla^g_\xi\xi=0,\ \nabla^g_\xi\phi=0$, so $(\Lie_\xi\phi)Y=(\nabla^g_\xi\phi)Y-\nabla^g_{\phi Y}\xi+\phi\nabla^g_Y\xi=-\bigl\{-\alpha\phi^2Y+\beta\phi Y\bigr\}+\phi\bigl\{-\alpha\phi Y+\beta(Y-\eta(Y)\xi)\bigr\}=0$, and hence  we have 
\begin{align}
\mathcal{L}_\xi\phi=0.\label{eq:L xi phi 0}
\end{align}

On a three-dimensional orientable Riemannian manifold $(M,g)$, a unit vector field $\xi$ defines an ACM structure as follows.
Set $\eta=\xi^\flat$ and we define an endomorphism field $\phi(X):=\xi\times X$.
Then we obtain $\phi^2=-\Id+\eta\otimes \xi$ and $g(\phi X,\phi Y)=g(X,Y)-\eta(X)\eta(Y)$.
Thus, $(\phi,\xi,\eta,g)$ is an ACM structure.

An ACM manifold $(M,\phi,\xi,\eta,g)$ is called \textcolor{blue}{\textbf{$\eta$-Einstein}} if there exist smooth real functions $\lambda,\nu$ such that $\Ric=\lambda g+\nu\eta\otimes\eta$.
In~\cite{Matsuno2026,Matsuno2025}, we classified simply connected three-dimensional trans-Sasakian $\eta$-Einstein Riemannian manifolds $(M,\phi,\xi,\eta,g)$ with complete $\xi$.
We review this classification in Section~\ref{sec:explicit-model-fields}.

\section{Rank dichotomy}\label{subsec:Unique continuation of rank-two}
In this section, we show that only two cases can occur for the rank of an RS field, namely, rank two almost everywhere or rank three almost everywhere.
This follows from the unique continuation of elliptic operators.

Let $(M,g)$ be a connected three-dimensional Riemannian spin manifold, and let $\Psi$ be an RS field.
Set $Z:=\{p\in M;\ \Psi(p)=0\}$ and $M^\circ:=M\setminus Z$.
Regard the spinor bundle as a real bundle $\Sigma_{\R}M$ of real rank $4$.
Let $\langle,\rangle$ be a spin invariant Hermitian inner product, and we can define a spin invariant real inner product $(,):=\Re\langle,\rangle$.
We can also define an orientation on $\Sigma_{\R}M$ induced from a complex structure of $\Sigma M$.
Let $*_{\Sigma}$ be the Hodge star operator determined by this metric and orientation.
Since the tangent bundle of an orientable three-dimensional manifold is trivial, a global orthonormal frame $\{e_1,e_2,e_3\}$ exists.
We define
\begin{equation}\label{cubic:normal-spinor}
 N:=*_{\Sigma}\bigl(\Psi(e_1)\wedge\Psi(e_2)\wedge\Psi(e_3)\bigr)
 \in\Gamma(\Sigma_{\R}M).
\end{equation}
The condition $N(p)=0$ is equivalent to $\rank_{\R}\Psi_p\le2$.
By Lemma~\ref{lem:rs-real-rank}, $\rank_{\R}\Psi_p\ne1$ thus the condition $N(p)=0$ is equivalent to $\rank_{\R}\Psi_p=2$ on $M^\circ$.
We have the following lemma.

\begin{lemma}\label{lem:N elliptic-inequality}
For every $p\in M^\circ$, there exist a neighborhood $U$ containing $p$ and a constant $C_U$ such that
\begin{equation}\label{cubic:elliptic-bound}
 \bigl|(\nabla^\Sigma)^*\nabla^\Sigma N\bigr|\le C_U\bigl(|N|+|\nabla^\Sigma N|\bigr)
\end{equation}
on $U$.
Here $(\nabla^\Sigma)^*\nabla^\Sigma$ is the rough Laplacian associated with the spin connection $\nabla^\Sigma$.
\end{lemma}
\begin{proof}
See Appendix~\ref{app:N elliptic-inequality}.
\end{proof}

\begin{proposition}\label{prop:rank-propagation}
Let $(M,g)$ be a connected three-dimensional Riemannian spin manifold, and let $\Psi$ be an RS field that is not identically zero.
Then either (R2) $\Psi$ has rank two almost everywhere on $M$, or (R3) $\Psi$ has rank three almost everywhere on $M$.
Furthermore, in case (R2), every point at which $\Psi$ does not have rank two is a zero of $\Psi$, and in case (R3), every point at which $\Psi$ does not have rank three is either a zero of $\Psi$ or a point at which $\Psi$ has rank two.
\end{proposition}
\begin{proof}
Since $D_{TM}$ is a Dirac-type operator and $\Psi\not\equiv0$, B\"ar's nodal set theorem~\cite{Bar1999} gives $\dim_{\mathcal H}Z\le1$.
Thus $Z$ has zero volume, and $M^\circ$ is connected by~\cite[Proposition~2.2]{Gerner2021}.
Lemma~\ref{lem:N elliptic-inequality} and~\cite[Theorem~2.1 and Remark~2.2]{BrandingEtAl} imply that $N$ satisfies the strong unique continuation property on $M^\circ$.

Suppose that $E:=\{N=0\}\cap M^\circ=\{p\in M:\operatorname{rank}_{\mathbb R}\Psi_p=2\}$ has positive volume, and choose a point $q\in E$ of density one.
On a neighborhood $U\subset M^\circ$ of $q$, choose coordinates $x$ with $x(q)=0$ and a smooth local frame $s_1,\ldots,s_4$ of $\Sigma_{\R}M$.
Write $N(p)=\sum_{\alpha=1}^4u^\alpha(x(p))s_\alpha(p)$ for $p\in U$, and set $u=(u^1,\ldots,u^4):x(U)\to\R^4$.
In these coordinates, $E\cap U$ is represented by $\{u=0\}$.
Since smooth coordinate changes preserve density-one points and the Riemannian volume density is smooth and positive, we have
\[
\lim_{r\downarrow0}\frac{\mathcal L^3(\{u=0\}\cap B_r(0))}{\mathcal L^3(B_r(0))}=1,
\]
where $\mathcal L^3$ denotes Lebesgue measure.
If there exists an integer $m\ge1$ such that $\partial^\beta u(0)=0\quad(\beta\in\mathbb N_0^3,\ |\beta|<m)$ and $P_m(x):=\sum_{|\beta|=m}\frac{\partial^\beta u(0)}{\beta!}x^\beta\not\equiv0$.
Then $\sup_{|x|\le1}|r^{-m}u(rx)-P_m(x)|=O(r)$ as $r\downarrow0$.
Hence there exists a constant $C>0$ such that, for every sufficiently small $r>0$, the estimate $|r^{-m}u(rx)-P_m(x)|\le Cr$ holds for all $x\in B_1(0)$.
At every point of this set, we have $|P_m(x)|=|P_m(x)-r^{-m}u(rx)|\le Cr$.
Consequently, we have
\begin{align*}
\frac{\mathcal L^3(\{u=0\}\cap B_r(0))}{\mathcal L^3(B_r(0))}
&=\frac{\mathcal L^3(\{x\in B_1(0):u(rx)=0\})}{\mathcal L^3(B_1(0))}\le\frac{\mathcal L^3(\{x\in B_1(0):|P_m(x)|\le Cr\})}{\mathcal L^3(B_1(0))}.
\end{align*}
Moreover, by continuity of measure for decreasing sequences of sets, we have $\lim_{r\downarrow0}\mathcal L^3(\{x\in B_1(0):|P_m(x)|\le Cr\})=\mathcal L^3(\{x\in B_1(0):P_m(x)=0\})=0$, because the zero set of the nonzero polynomial $P_m$ has Lebesgue measure zero.
This contradicts the density-one property of $q$.
Thus $N$ vanishes to infinite order at $q$, and strong unique continuation gives $N\equiv0$ on $M^\circ$.

If $N\equiv0$ on $M^\circ$, Lemma~\ref{lem:rs-real-rank} gives $\rank_{\R}\Psi=2$ throughout $M^\circ$, and hence $\{p\in M:\rank_{\R}\Psi_p\ne2\}=Z$.
Otherwise $E$ has zero volume, so $\rank_{\R}\Psi=3$ almost everywhere on $M^\circ$, and Lemma~\ref{lem:rs-real-rank} gives $\{p\in M:\rank_{\R}\Psi_p\ne3\}=Z\cup E$.
\end{proof}

\section{Geometry determined by a rank-two Rarita-Schwinger field}
\label{sec:kernel-geometry}
In this section, we prove one of the main results, Theorem~\ref{thm:removable-main}.
In Subsection~\ref{subsec:A trans-Sasakian eta-Einstein structure on a rank-two domain}, we show that on a rank-two domain, the line bundle $\ker\Psi$ defines a trans-Sasakian $\eta$-Einstein structure.
In this case, $\ker\Psi$ determines the characteristic vector field of the trans-Sasakian structure up to sign.
At this stage, the trans-Sasakian $\eta$-Einstein structure is defined away from the zeros of $\Psi$, but at the zeros, $\ker\Psi$ does not specify a line, so the characteristic vector field is not determined by $\Psi$.
In Subsection~\ref{subsec:Extension of trans-Sasakian eta-Einstein structure onto a rank-zero domain}, we show that the trans-Sasakian $\eta$-Einstein structure extends uniquely to the zeros of $\Psi$.

\subsection{A trans-Sasakian \texorpdfstring{$\eta$}{eta}-Einstein structure on a rank-two domain}\label{subsec:A trans-Sasakian eta-Einstein structure on a rank-two domain}
In this subsection, we show that on a domain where an RS field $\Psi$ has rank two, a trans-Sasakian $\eta$-Einstein structure is determined up to sign.

\begin{lemma}\label{lem:kernel-form}
On a three-dimensional Riemannian spin manifold $(M,g)$, assume that $\Psi\in\Gamma(\Sigma_{3/2}M)$ has rank two on an open set $U$.
Let $\xi$ be a vector field on $U$ satisfying $\ker_\R\Psi=\R \xi,\ ||\xi||^2=1$, and let $\{e_1=\xi,e_2,e_3\}$ be an orthonormal frame on $U$.
Then $\Phi:=\Psi(e_2)\ne0$ and $\Psi(e_3)=\gamma_1\Phi$, and therefore we have $\Psi=\Phi\otimes\theta^2+\gamma_1\Phi\otimes\theta^3$, where $\{\eta,\theta^2,\theta^3\}$ is the corresponding coframe on $U$.
\end{lemma}
\begin{proof}
Since $e_1=\xi\in\ker\Psi$, we have $0=\Pi(\Psi) =\gamma_1\Psi(e_1)+\gamma_2\Psi(e_2)+\gamma_3\Psi(e_3)=\gamma_2\Phi+\gamma_3\Psi(e_3)$, which implies $\Psi(e_3)=-\gamma_2\gamma_3\Phi=\gamma_1\Phi$.
If $\Phi(p)=0$ at some point $p\in U$, then $\Psi_p(e_1(p))=\Psi_p(e_2(p))=\Psi_p(e_3(p))=0$, so $\rank_\R\Psi_p=0$, contradicting $\rank_\R\Psi_p=2$.
Hence $\Phi$ is nowhere zero.
\end{proof}

\begin{proposition}\label{prop:rank-two-kernel-geometry}
On a three-dimensional Riemannian spin manifold $(M,g)$, let $\Psi\in\Gamma(\Sigma_{3/2}M)$ be an RS field that has rank two on an open set $U$.
Let $\xi$ be a smooth unit vector field on $U$ satisfying $\ker_\R\Psi=\R \xi,\ ||\xi||^2=1$.
Then the ACM structure defined by $\xi$ on $U$ is trans-Sasakian and $\eta$-Einstein.
\end{proposition}
\begin{proof}
Let $\{e_1=\xi,e_2,e_3\}$ be an orthonormal frame on $U$, and put $a_i=g(\nabla^g_{e_i}\xi,e_2),\allowbreak\ b_i=g(\nabla^g_{e_i}\xi,e_3)$.
Since we have \(\Phi=\Psi(e_2)\ne0,\ \Psi(\xi)=0,\ \Psi(e_3)=\gamma_1\Phi\), we have $0=(D_{TM}\Psi)(\xi)=\sum_i\gamma_i\left\{\nabla^\Sigma_{e_i}(\Psi(\xi))-\Psi(\nabla^g_{e_i}\xi)\right\}=-\sum_i\gamma_i(a_i\Phi+b_i\gamma_1\Phi)=-(-b_1\operatorname{Id}+a_1\gamma_1 +(a_2-b_3)\gamma_2+(a_3+b_2)\gamma_3)\Phi$.
Thus we obtain $a_1=b_1=0,\ a_2=b_3=:\beta,\ a_3=-b_2=:\alpha$ from the real linear independence of \(\Phi,\gamma_1\Phi,\gamma_2\Phi,\gamma_3\Phi\).
Therefore $\nabla^g_\xi\xi=0,\ \nabla^g_{e_2}\xi=\beta e_2-\alpha e_3,\ \nabla^g_{e_3}\xi=\alpha e_2+\beta e_3$ holds, which implies \eqref{eq:xi-connection} holds, then the ACM structure defined by $\xi$ is trans-Sasakian.

Next, we show that the structure is $\eta$-Einstein.
Setting $r_{ij}=\Ric(e_i,e_j)$, we obtain
\begin{align*}
 0
 &=\cR(\Psi)
 =\sum_{i=1}^3\gamma(\Ric(e_i))\Psi(e_i)
 =\gamma(\Ric(e_2))\Phi
   +\gamma(\Ric(e_3))\gamma_1\Phi\\
 &=\bigl[r_{21}\gamma_1+r_{22}\gamma_2+r_{23}\gamma_3
   +(r_{31}\gamma_1+r_{32}\gamma_2+r_{33}\gamma_3)\gamma_1\bigr]\Phi\\
 &=\bigl[r_{21}\gamma_1+r_{22}\gamma_2+r_{23}\gamma_3
   +r_{31}\gamma_1^2+r_{32}\gamma_2\gamma_1
   +r_{33}\gamma_3\gamma_1\bigr]\Phi\\
 &=\bigl[r_{21}\gamma_1+r_{22}\gamma_2+r_{23}\gamma_3
   -r_{31}\Id+r_{32}\gamma_3-r_{33}\gamma_2\bigr]\Phi\\
 &=\bigl[-r_{31}\Id+r_{21}\gamma_1
   +(r_{22}-r_{33})\gamma_2+(r_{23}+r_{32})\gamma_3\bigr]\Phi\\
 &=\bigl[-r_{31}\Id+r_{21}\gamma_1
   +(r_{22}-r_{33})\gamma_2+2r_{23}\gamma_3\bigr]\Phi.
\end{align*}
The real linear independence gives $r_{21}=r_{31}=r_{23}=0,\ r_{22}=r_{33}$, so the structure is $\eta$-Einstein.
\end{proof}

\begin{proposition}\label{prop:rank-two-local-geometry}
Let $(M,g)$ be a connected three-dimensional Riemannian spin manifold, and let $\Psi\in\Gamma(\Sigma_{3/2}M)$ be an RS field of type (R2).
Let $M^\circ=M\setminus\{\Psi=0\}$, and $L^\circ:=\ker_\R\Psi|_{M^\circ}$ is a smooth real line bundle.
Then every local unit section $\xi$ of $L^\circ$ defines a trans-Sasakian $\eta$-Einstein structure on its domain.
\end{proposition}
\begin{proof}
By Proposition~\ref{prop:rank-propagation}, $\rank_\R\Psi=2$ on $M^\circ$.
Thus, $L^\circ$ is a smooth real line bundle and admits a smooth unit section on a neighborhood of each point.
Applying Proposition~\ref{prop:rank-two-kernel-geometry} to such a local unit section $\xi$ gives the conclusion.
\end{proof}

\subsection{Extension of the trans-Sasakian \texorpdfstring{$\eta$}{eta}-Einstein structure across the zero set}\label{subsec:Extension of trans-Sasakian eta-Einstein structure onto a rank-zero domain}
In this subsection, we extend the one-dimensional distribution $L^\circ$ on $M^\circ$ to all of $M$ and complete the proof of Theorem~\ref{thm:removable-main}.

\begin{lemma}\label{lem:rs-norm-transport}
Let $(M,g)$ be a connected three-dimensional Riemannian spin manifold, and let $\Psi$ be an RS field of type (R2).
Let $M^\circ=M\setminus\{\Psi=0\}$ and $\xi$ be a local unit section of the one-dimensional distribution $L^\circ:=\ker_\R\Psi|_{M^\circ}$ on $M^\circ$, and let $(\alpha,\beta)$ be the type of the trans-Sasakian structure defined by $\xi$.
Then $f:=|\Psi|$ satisfies $\xi\log f=-\beta$ on the domain of $\xi$.
Moreover, $Z:=f^{-1}(0)$ satisfies $\mathcal H^1(Z\cap K)<\infty$ for every compact set $K$, where $\mathcal H^1$ is one-dimensional Hausdorff measure.
\end{lemma}
\begin{proof}
Choose a local orthonormal frame $\{\xi,e_2,e_3\}$ on $M^\circ$, with dual coframe $\{\eta,\theta^2,\theta^3\}$.
Setting $\Phi:=\Psi(e_2)$, we can write $\Psi=\Phi\otimes\theta^2+\gamma_1\Phi\otimes\theta^3$.
Thus, $||\Psi||^2=2||\Phi||^2$.
Although $\Phi$ is in general only locally defined, $|\Phi|=\frac{|\Psi|}{\sqrt2}$ is defined on all of $M$.
Write the connection coefficients as $a_i=g(\nabla^g_{e_i}\xi,e_2),\ b_i=g(\nabla^g_{e_i}\xi,e_3),\ c_i=g(\nabla^g_{e_i}e_2,e_3)$.
Since $[\xi,e_2]=\nabla^g_\xi e_2-\nabla^g_{e_2}\xi=-a_1\xi-a_2e_2+(c_1-b_2)e_3$, we have
$0=(\dd^\nabla\Psi)(\xi,e_2)=\nabla^\Sigma_\xi(\Psi(e_2))-\nabla^\Sigma_{e_2}(\Psi(\xi))-\Psi([\xi,e_2])
=\nabla^\Sigma_\xi\Phi-\Psi\bigl(-a_1\xi-a_2e_2+(c_1-b_2)e_3\bigr)
=\nabla^\Sigma_\xi\Phi+a_2\Phi-(c_1-b_2)\gamma_1\Phi$, thus we obtain $\nabla^\Sigma_\xi\Phi=-a_2\Phi+(c_1-b_2)\gamma_1\Phi=-\beta\Phi+(c_1+\alpha)\gamma(\xi)\Phi$, so $\xi|\Phi|^2=2\Re\langle\nabla^\Sigma_\xi\Phi,\Phi\rangle =-2\beta|\Phi|^2$, and hence $\xi\log |\Phi|=-\beta$.
Since $D_{TM}$ is a Dirac-type operator, B\"ar's theorem on the zero set of a nontrivial solution~\cite{Bar1999} implies that $Z:=f^{-1}(0)=\{\Psi=0\}$ satisfies $\mathcal H^1(Z\cap K)<\infty$ for every compact set $K$.
\end{proof}

Since the set $Z$ on which $L^\circ$ is not defined is a ``small'' set, one expects $L^\circ$ to extend across $Z$ as in a removable singularity problem.
Indeed, the following lemma shows that $L^\circ$ extends across $Z$.

\begin{lemma}\label{lem:geom}
Let $Z\subset M$ be a closed set such that $\mathcal H^1(Z\cap K)<\infty$ for every compact set $K$, and let $M^\circ=M\setminus Z$.
Suppose that there exists a continuous function $f:M\to[0,\infty)$ satisfying $f^{-1}(0)=Z$ and $f\in C^\infty(M^\circ)$, and that $M^\circ$ carries a smooth line distribution $L^\circ$ whose local unit sections $\xi$ define trans-Sasakian $\eta$-Einstein structures.
Also assume that a local unit section $\xi$ satisfies $\xi(\log f)=-\beta$ on a domain of $\xi$, where $(\alpha,\beta)$ is the type of the trans-Sasakian structure defined by $\xi$.
Then $L^\circ$ extends uniquely to a smooth real line bundle $L$ on all of $M$.
Every local unit section $\xi$ of the extended bundle $L$ defines a trans-Sasakian $\eta$-Einstein structure on its domain.
\end{lemma}
\begin{proof}
See Appendix~\ref{app:geom}.
\end{proof}

The first main result now follows immediately from these results.
\begin{theorem}\label{thm:removable-main}
Let $(M,g)$ be a connected and simply connected three-dimensional Riemannian spin manifold, and let $\Psi$ be an RS field of type (R2), that is, it has rank two almost everywhere on $M$.
Then $M$ admits a trans-Sasakian $\eta$-Einstein structure $(\phi,\xi,\eta,g)$ such that $\ker_\R\Psi_p=\R\xi_p$ at every point $p$ where $\Psi(p)\ne0$.
The vector field $\xi$ satisfying this condition is unique up to sign.
\end{theorem}
\begin{proof}
Set $Z:=\{\Psi=0\}$ and $M^\circ:=M\setminus Z$.
By Proposition~\ref{prop:rank-two-local-geometry}, $\rank_\R\Psi=2$ on $M^\circ$, and every local unit section of $L^\circ:=\ker_\R\Psi|_{M^\circ}$ defines a trans-Sasakian $\eta$-Einstein structure.
By Lemma~\ref{lem:rs-norm-transport}, $f:=|\Psi|$ and $Z$ satisfy the assumptions of Lemma~\ref{lem:geom}.
Thus, by Lemma~\ref{lem:geom}, $L^\circ$ extends uniquely to a smooth real line bundle $L\subset TM$ on $M$, and each of its local unit sections defines a trans-Sasakian $\eta$-Einstein structure.

Let $S(L)$ be a unit sphere bundle of $L$ and $S(L)\to M$ be a double cover.
Since $M$ is simply connected, $S(L)\simeq M\times\{\pm1\}$.
Hence a global smooth unit section $\xi$ exists.
\end{proof}

\section{Rarita-Schwinger fields on trans-Sasakian $\eta$-Einstein manifolds}\label{sec:Rarita-Schwinger fields on trans-Sasakian eta-Einstein manifolds}

\subsection{Kernel rigidity}
In this subsection, we show that the kernel of an RS field on a three-dimensional trans-Sasakian $\eta$-Einstein manifold whose metric does not have constant curvature contains the characteristic vector field.

\begin{lemma}\label{lem:kernel-ucp}
Let $(M,\phi,\xi,\eta,g)$ be a three-dimensional trans-Sasakian manifold of type $(\alpha,\beta)$, and let $\Psi$ be an RS field.
Set $\chi:=\Psi(\xi)$.
Then, on every relatively compact open set $U$, there is a constant $C_U$ such that
\begin{equation}\label{eq:chi-bound}
|(\nabla^\Sigma)^*\nabla^\Sigma\chi|\le C_U(|\nabla^\Sigma\chi|+|\chi|).
\end{equation}
\end{lemma}
\begin{proof}
Write $\psi_j=\Psi(e_j)$ and $e_1=\xi$.
Substituting $D_{TM}\Psi=0$ and the identities $\nabla^g_\xi\xi=0$, $\nabla^g_{e_2}\xi=\beta e_2-\alpha e_3$, and $\nabla^g_{e_3}\xi=\beta e_3+\alpha e_2$ from~\eqref{eq:xi-connection} into $D_\Sigma\chi=(D_{TM}\Psi)(\xi)+\sum_i\gamma_i\Psi(\nabla^g_{e_i}\xi)$, we obtain $D_\Sigma\chi=-\beta\gamma_1\chi   -\alpha(\gamma_2\psi_3-\gamma_3\psi_2)$.
Since $\Pi(\Psi)=0$ gives $\gamma_2\psi_3-\gamma_3\psi_2=\chi$, we obtain $(D_\Sigma+\beta\gamma(\xi)+\alpha)\chi=0$.
Set $A:=\beta\gamma(\xi)+\alpha\Id$ and apply $D_\Sigma$ once more to $D_\Sigma\chi=-A\chi$.
The Lichnerowicz formula $D_\Sigma^2=(\nabla^\Sigma)^*\nabla^\Sigma+\Scal/4$ and the identity $D_\Sigma(A\chi)=\sum_i\gamma_i\{(\nabla^g_{e_i}A)\chi+A\nabla^\Sigma_{e_i}\chi\}$ give
\begin{align}
 (\nabla^\Sigma)^*\nabla^\Sigma\chi
 =-\frac{\Scal}{4}\chi-D_\Sigma(A\chi)
 =-\frac{\Scal}{4}\chi
   -\sum_i\gamma_i(\nabla^g_{e_i}A)\chi
   -\sum_i\gamma_i A\nabla^\Sigma_{e_i}\chi.
\end{align}
Here $\nabla^g_XA=X(\beta)\gamma(\xi)+\beta\gamma(\nabla^g_X\xi)+X(\alpha)\Id$.
On any relatively compact $U$, the quantities $\Scal,A,\nabla^g A$ are bounded, and the estimates $|\gamma(e_i)\zeta|=|\zeta|$ and $\sum_i|\nabla^\Sigma_{e_i}\chi|\le\sqrt3|\nabla^\Sigma\chi|$ give
\begin{align*}
 |(\nabla^\Sigma)^*\nabla^\Sigma\chi|
 \le \left(\frac{|\Scal|}{4}+\sum_i\|\nabla^g_{e_i}A\|\right)|\chi|
   +\sqrt3\|A\|\,|\nabla^\Sigma\chi|
 \le C_U(|\chi|+|\nabla^\Sigma\chi|).
\end{align*}
\end{proof}

The following is the second main result of this paper.

\begin{theorem}\label{thm:fixed}
Let $(M,\phi,\xi,\eta,g)$ be a connected three-dimensional trans-Sasakian $\eta$-Einstein manifold which does not have constant curvature.  
Every RS field $\Psi$ that is not identically zero is of type (R2) and satisfies $\ker_{\R}\Psi_p=\R\xi_p$ for $\Psi_p\ne0$.
\end{theorem}
\begin{proof}
Set $\chi:=\Psi(\xi)$ and write $\Ric=ag+b\eta\otimes\eta$.
By Lemma~\ref{lem:ricci} and $\Pi\Psi=0$, we have $ 0=a\sum_i\gamma_i\Psi(e_i)+b\sum_i\eta(e_i)\gamma(\xi)\Psi(e_i)=a\,\Pi\Psi+b\gamma(\xi)\Psi(\xi)=b\gamma(\xi)\chi$.
Since the manifold does not have constant curvature, $b\ne0$ on some nonempty open set $U$.
Thus $\chi=0$ on $U$.
By Lemma~\ref{lem:kernel-ucp}, there is locally a constant $C>0$ such that
\[
|-g^{ij}\partial_i\partial_j\chi_a|\le C\bigl(|\chi_1|+|\chi_2|+|\partial\chi_1|+|\partial\chi_2|\bigr),\qquad a=1,2.
\]
The vector-valued version of Aronszajn's unique continuation theorem~\cite[Theorem~2.1]{BrandingEtAl} and the connectedness of $M$ therefore imply that $\chi\equiv0$ on all of $M$.
Consequently, $\rank_\R\Psi_p\le2$, and Lemma~\ref{lem:rs-real-rank} gives $\rank_\R\Psi_p\in\{0,2\}$, with $\ker_\R\Psi_p=\R\xi_p$ at each point $p$ such that $\Psi_p\ne0$.
\end{proof}

\subsection{Local existence and global characterization}
In this subsection, we prove local existence of RS fields of type (R2) satisfying $\Psi(\xi)=0$ and determine the space of all global RS fields ot type (R2).
Let $(M^3,\phi,\xi,\eta,g)$ be a three-dimensional trans-Sasakian manifold.
If the manifold is $\eta$-Einstein and its metric does not have constant curvature, the results of the previous subsection imply that every RS field satisfies $\Psi(\xi)=0$.
We define the subbundle $\mathcal R_\xi\subset\Sigma_{3/2}M$ as follows.
\begin{equation}
 \mathcal R_\xi:=\left\{A\in \Sigma_{3/2}M\ ;\ A(\xi)=0\right\}
\end{equation}
Define the eigenline bundles of the spinor bundle for $\gamma(\xi)$ by $ \Sigma_\pm:=\ker\bigl(\gamma(\xi)\mp i\Id\bigr)$, and put
\begin{align*}
\mathcal{R}_\xi^\pm&:=\mathcal{R}_\xi\cap(T^*M\otimes\Sigma_\pm)
=\left\{A\in \mathcal{R}_\xi\ ;\ 
 \gamma(\xi)A(X)=\pm iA(X)\ \text{for all }X\in TM\right\}
\end{align*}
Then $\mathcal R_\xi=\mathcal R_\xi^+\oplus\mathcal R_\xi^-$.
Let $\mathcal H:=(\ker\eta)\otimes\C\subset TM^\C$, and write its eigenspace decomposition for $\phi$ as $\mathcal H=\mathcal H^{1,0}\oplus\mathcal H^{0,1}$.
Here $\mathcal H^{1,0}$ is the $+i$ eigenspace, and $\mathcal H^{0,1}$ is the $-i$ eigenspace.
Define the subbundle $\mathcal V$ of $TM^\C$ by $\mathcal V:=\C\xi\oplus\mathcal H^{0,1}$.
Its annihilator is $\operatorname{Ann}(\mathcal V)=(\mathcal{H}^{1,0})^*$.
We denote $\mathcal{V}_+=\mathcal{V},\ \mathcal{V}_-=\overline{\mathcal{V}}$.
The following holds for $\mathcal{V}_\pm$ and $\mathcal{R}_\xi^\pm$.

\begin{lemma}
Let $(M^3,\phi,\xi,\eta,g)$ be a three-dimensional trans-Sasakian manifold.
The distribution $\mathcal V_\pm$ is involutive.
For $Z\in\mathcal H^{0,1}, \sigma\in\Sigma_+$ or $Z\in\mathcal H^{1,0}, \sigma\in\Sigma_-$, we have
\begin{align}
\gamma(Z)\sigma=0,\label{eq:Zsigma0}
\end{align}
and furthermore $\mathcal R_\xi^+=\Sigma_+\otimes (\mathcal{H}^{1,0})^*$ and $\mathcal R_\xi^-=\Sigma_-\otimes(\mathcal{H}^{0,1})^*$ hold.
\end{lemma}
\begin{proof}
For \(X\in\Gamma(\mathcal H)\), we have $\eta([\xi,X])=-g(X,\nabla^g_\xi\xi)=0$, so $[\xi,X]\in\Gamma(\mathcal H)$.
By \eqref{eq:L xi phi 0}, $[\xi,X+i\phi X]=[\xi,X]+i[\xi,\phi X]=[\xi,X]+i\phi[\xi,X]\in\Gamma(\mathcal H^{0,1})$ holds, and hence $\mathcal{V}_+$ is involutive.

For an orthonormal frame $\{e_1=\xi,e_2,e_3=\phi e_2\}$, we can write $Z\in\mathcal H^{0,1}$ as $Z=z(e_2+i e_3),\ z\in C^\infty(M,\C)$, and thus $\gamma(Z)\sigma=z(\gamma_2+i\gamma_3)\sigma=z\gamma_2(\operatorname{Id}+i\gamma_1)\sigma=z\gamma_2(\sigma-\sigma)=0$ holds.
Also, for any $0\ne W\in\mathcal H^{1,0}$ and $0\ne\sigma\in\Sigma_+$, we have $\gamma(W)\sigma\ne0$.
Thus $\ker\bigl(\mathcal H\ni Z\longmapsto\gamma(Z)\sigma\bigr)=\mathcal H^{0,1}$.
Any $A\in\mathcal{R}_\xi^+$ can be written, using $\sigma\in\Gamma(\Sigma_+)$, as $A=\sigma\otimes\lambda$.
Here $\lambda\in T^*M^\C$ satisfies $\lambda(\xi)=0$.
Then $\Pi(A)=0\ \Leftrightarrow\ \gamma(\lambda^\sharp)\sigma=0\ \Leftrightarrow\ \lambda^\sharp\in\mathcal H^{0,1}\ \Leftrightarrow\ \lambda\in (\mathcal{H}^{1,0})^*$.
Thus $\mathcal R_\xi^+=\Sigma_+\otimes (\mathcal{H}^{1,0})^*$.
The assertion for $\mathcal{R}_\xi^-$ follows similarly.
\end{proof}

Since $\mathcal V_+\wedge(TM\otimes\C)$ generates $\Lambda^2(TM\otimes\C)$, for $\Psi_+\in\Gamma(\mathcal R_\xi^+)$ the condition $\dd^\nabla\Psi_+=0$ is equivalent to $(\dd^\nabla\Psi_+)(X,Y)=0$ for all $X\in\Gamma(\mathcal V_+)$ and $Y\in\Gamma(TM^\C)$.
Furthermore, since $\mathcal R_\xi^+=\Sigma_+\otimes(\mathcal H^{1,0})^*$, we have $\Psi_+(X)=0$ for $X\in\Gamma(\mathcal V_+)$, and hence
\begin{align*}
(\dd^\nabla\Psi_+)(X,Y)=\nabla^\Sigma_X(\Psi_+(Y))-\Psi_+([X,Y]).
\end{align*}
By Lemma~\ref{lem:closed}, $\Psi_+$ is an RS field if and only if the right-hand side vanishes for all such $X$ and $Y$.
We interpret this condition as parallelism with respect to a partial connection on $\mathcal R_\xi^+$.

We follow~\cite[Section~3.2]{BiswasSchumacher} for the definition of a \textcolor{blue}{\textbf{partial connection}}.
A partial connection on a complex vector bundle $\mathcal W$ along a tangent distribution $\mathcal{V}$ is a connection defined by a $\C$-linear operator $\mathcal D:\Gamma(\mathcal W)\longrightarrow\Gamma(\mathcal{V}^*\otimes\mathcal W)$ extended by the Leibniz rule $\mathcal D(fs)=\dd f|_{\mathcal{V}_+}\otimes s+f\mathcal Ds$, allowing differentiation only in the $\mathcal{V}$ directions.
Its curvature is defined by $F_{\mathcal D}(X,Y)s:=\mathcal D_X\mathcal D_Ys-\mathcal D_Y\mathcal D_Xs-\mathcal D_{[X,Y]}s,\ (X,Y\in\Gamma(\mathcal{V}))$.

The \textcolor{blue}{\textbf{Bott connection}}~\cite[Definition~3.5]{Hawkins2008} is a partial connection along an involutive tangent distribution.
We define the Bott connection as a partial connection on $(\mathcal{H}^{1,0})^*$ along $\mathcal{V}_+$ by the following formula.
\begin{equation}\label{eq:bott-connection}
 \nabla^B_X\lambda:=\Lie_X\lambda,\quad X\in\Gamma(\mathcal{V}_+),\ \lambda\in\Gamma((\mathcal{H}^{1,0})^*)
\end{equation}
Although it is not obvious that this defines a partial connection, the involutivity of $\mathcal{V}_+$ guarantees it.
Indeed, since $\lambda(X)=0$, we have $\Lie_{fX}\lambda=\iota_{fX}\dd\lambda=f\Lie_X\lambda$.
Also, for $X,Y\in\Gamma(\mathcal{V}_+)$, we have $(\mathcal{L}_X\lambda)(Y)=X(\lambda(Y))-\lambda([X,Y])=0$, so $\mathcal{L}_X\lambda\in\Gamma((\mathcal{H}^{1,0})^*)$.
Since it is defined by Lie derivatives, we have $F_{\nabla^B}(X,Y)\lambda =(\Lie_X\Lie_Y-\Lie_Y\Lie_X-\Lie_{[X,Y]})\lambda=0$, so the Bott connection is flat.

We can also see that the spin connection preserves $\Sigma_+$ along $\mathcal{V}_+$ as follows.
Writing \(X\in\mathcal{V}_+\) as \(X=a\xi+Z\), with \(Z\in\mathcal H^{0,1}\), equation~\eqref{eq:xi-connection} gives $\nabla^g_X\xi=-\alpha\phi Z+\beta Z=(\beta+i\alpha)Z$, so $\nabla^g_X\xi\in\mathcal H^{0,1}$.
For $\sigma\in\Gamma(\Sigma_+)$, we have $0=\nabla^\Sigma_X\bigl(\gamma(\xi)\sigma-i\sigma\bigr)=\gamma(\nabla^g_X\xi)\sigma+\gamma(\xi)\nabla^\Sigma_X\sigma-i\nabla^\Sigma_X\sigma$, so $(\gamma(\xi)-i)\nabla^\Sigma_X\sigma=-\gamma(\nabla^g_X\xi)\sigma=0$ by~\eqref{eq:Zsigma0}, and hence the spin connection preserves $\Sigma_+$ along $\mathcal{V}_+$.
Similarly, the spin connection preserves $\Sigma_-$ along $\mathcal{V}_-$.

Using the tensor product of the spin connection and the Bott connection, we define a partial connection on $\mathcal R_\xi^+$ along $\mathcal{V}_+$ by
\begin{align*}
\mathcal D^+_X(\sigma\otimes\lambda):=(\nabla^\Sigma_X\sigma)\otimes\lambda +\sigma\otimes\nabla^B_X\lambda,\\
\sigma\in\Gamma(\Sigma_+),\ \lambda\in\Gamma((\mathcal{H}^{1,0})^*),\ X\in\Gamma(\mathcal{V}_+)
\end{align*}
Similarly, we define a partial connection on $\mathcal R_\xi^-$ along $\mathcal{V}_-$, denoted by $\mathcal D^-$.
By this definition, for $\Psi_\pm=\sigma\otimes\lambda\in\mathcal{R}_\xi^\pm$, we obtain 
\begin{align*}
(\mathcal D_X^\pm\Psi_\pm)(Y)
&=\bigl((\nabla^\Sigma_X\sigma)\otimes\lambda
+\sigma\otimes\nabla^B_X\lambda\bigr)(Y)
=\lambda(Y)\nabla^\Sigma_X\sigma
+(\mathcal L_X\lambda)(Y)\sigma\\
&=\lambda(Y)\nabla^\Sigma_X\sigma
+\bigl\{X(\lambda(Y))-\lambda([X,Y])\bigr\}\sigma
=\nabla^\Sigma_X\bigl(\lambda(Y)\sigma\bigr)
-\lambda([X,Y])\sigma\\
&=\nabla^\Sigma_X(\Psi_\pm(Y))-\Psi_\pm([X,Y])=\dd^\nabla\Psi_\pm(X,Y)
\end{align*}
Then we have the following characterization:

\begin{proposition}\label{prop:intrinsic-rs-connection}
A section $\Psi:=\Psi_++\Psi_-\in\Gamma(\mathcal R_\xi),\ \Psi_\pm\in\Gamma(\mathcal R_\xi^\pm)$ is an RS field if and only if $ \mathcal D^\pm_X\Psi_\pm=0,\ X\in\Gamma(\mathcal{V}_\pm)$.
\end{proposition}
\begin{proof}
Since $(\mathcal D^\pm\Psi_\pm)=0\ \Leftrightarrow\ \dd^\nabla\Psi_\pm=0$ and $\Sigma_+\cap\Sigma_-=\{0\}$, the condition $\dd^\nabla\Psi=0$ is equivalent to $\mathcal{D}^\pm\Psi_\pm=0$.
\end{proof}

The relation between $\mathcal{R}_\xi^\pm$ and $\mathcal{D}^\pm$ is given as follows.

\begin{lemma}\label{lemma:conjugate}
We have $\overline{\mathcal R_\xi^+}\simeq\mathcal R_\xi^-$, and under this identification, $\overline{\mathcal D^+}=\mathcal D^-$.
\end{lemma}
\begin{proof}
The parallel quaternionic structure $\mathfrak j:\Sigma M\to\Sigma M$ of the three-dimensional spin representation is antilinear, satisfies $\mathfrak j^2=-\Id$, and commutes with Clifford multiplication by real vectors, and $\mathfrak j$ exchanges $\Sigma_+$ and $\Sigma_-$.
Thus we obtain the following complex linear bundle isomorphism.
\[
 \begin{aligned}
 I:\overline{\mathcal R_\xi^+}
 =\overline{\Sigma_+}\otimes
       \overline{(\mathcal H^{1,0})^*}
 &\longrightarrow
 \Sigma_-\otimes(\mathcal H^{0,1})^*
 =\mathcal R_\xi^-,\\
 I(\overline\sigma\otimes\overline\lambda)
 &:=\mathfrak j\sigma\otimes\overline\lambda
 \end{aligned}
\]
The conjugate connection $\overline{\mathcal D^+}$ is defined as a partial connection on $\overline{\mathcal R_\xi^+}$ along $\mathcal{V}_-=\overline{\mathcal{V}_+}$ by $(\overline{\mathcal D^+})_{\overline X}\overline s :=\overline{\mathcal D^+_Xs},\ X\in\Gamma(\mathcal{V}_+),\  s\in\Gamma(\mathcal R_\xi^+)$.
Then the following equality holds: $\mathcal D^-_{\overline X}I(\overline\sigma\otimes\overline\lambda)=(\nabla^\Sigma_{\overline X}\mathfrak j\sigma)\otimes\overline\lambda +\mathfrak j\sigma\otimes\Lie_{\overline X}\overline\lambda =\mathfrak j(\nabla^\Sigma_X\sigma)\otimes\overline\lambda  +\mathfrak j\sigma\otimes\overline{\Lie_X\lambda}=I\bigl(\overline{\mathcal D^+_X(\sigma\otimes\lambda)}\bigr) =I\bigl((\overline{\mathcal D^+})_{\overline X}(\overline\sigma\otimes\overline\lambda)\bigr)$.
\end{proof}

On three-dimensional trans-Sasakian $\eta$-Einstein manifolds, we may investigate a parallel section of $\mathcal{R}^\pm_\xi$ with respect to $\mathcal{D}^\pm$ to find an RS field.
The key fact is as follows.

\begin{proposition}\label{prop:intrinsic-rs-flatness}
If $(M,\phi,\xi,\eta,g)$ is $\eta$-Einstein, then $\mathcal D^\pm$ are flat.
\end{proposition}
\begin{proof}
We show that the $\eta$-Einstein condition implies that the spin curvature restricted to $\mathcal V_+\times\mathcal V_+$ vanishes on $\Sigma_+$.
The well-known curvature contraction identity $\gamma(\Ric(X))=2\sum_{i=1}^3\gamma_i\mathrm{Riem}(e_i,X)$ gives $\Scal_g\Id=-2\sum_{i,j=1}^3\gamma_i\gamma_j\mathrm{Riem}(e_j,e_i)$, for instance \cite{Friedrich-Kim}.
Setting $T:=\gamma_1\mathrm{Riem}(e_2,e_3)+\gamma_2\mathrm{Riem}(e_3,e_1)+\gamma_3\mathrm{Riem}(e_1,e_2)$, we obtain
$\Scal_g\Id=-2\sum_{i,j=1}^3\gamma_i\gamma_j\mathrm{Riem}(e_j,e_i)=4\bigl(\gamma_2\gamma_3\mathrm{Riem}(e_2,e_3)+\gamma_3\gamma_1\mathrm{Riem}(e_3,e_1)+\gamma_1\gamma_2\mathrm{Riem}(e_1,e_2)\bigr)=-4\bigl(\gamma_1\mathrm{Riem}(e_2,e_3)+\gamma_2\mathrm{Riem}(e_3,e_1)+\gamma_3\mathrm{Riem}(e_1,e_2)\bigr)=-4T$.
Applying the same contraction identity with $X=e_1$ gives
$\frac12\gamma_1\gamma(\Ric(e_1))=\gamma_1\sum_{i=1}^3\gamma_i\mathrm{Riem}(e_i,e_1)=\gamma_1(-\gamma_2\mathrm{Riem}(e_1,e_2)+\gamma_3\mathrm{Riem}(e_3,e_1))=\gamma_3\mathrm{Riem}(e_1,e_2)+\gamma_2\mathrm{Riem}(e_3,e_1)=T-\gamma_1\mathrm{Riem}(e_2,e_3)=-\frac{\Scal_g}{4}\Id-\gamma_1\mathrm{Riem}(e_2,e_3)$.
Hence we obtain $\mathrm{Riem}(e_2,e_3)=\frac{\Scal_g}{4}\gamma_1-\frac12\gamma(\Ric(e_1)),\ \mathrm{Riem}(e_3,e_1)=\frac{\Scal_g}{4}\gamma_2-\frac12\gamma(\Ric(e_2))$ and $\mathrm{Riem}(e_1,e_2)=\frac{\Scal_g}{4}\gamma_3-\frac12\gamma(\Ric(e_3))$.
These equations can be written in the following form.
\[
 \mathrm{Riem}(X,Y)
 =\frac{\Scal_g}{4}\gamma(X\times Y)
       -\frac12\gamma\bigl(\Ric(X\times Y)\bigr)
\]
Since $\Ric=\lambda g+\nu\eta\otimes\eta$, for $Z\in\mathcal H^{0,1}$ we have $\Ric^\sharp(Z)=\lambda Z$.
Noting that $\xi\times Z=\phi Z=-iZ$, for $\sigma\in\Sigma_+$ we have $\mathrm{Riem}(\xi,Z)\sigma=-i\left(\frac{\Scal_g}{4}-\frac{\lambda}{2}\right)\gamma(Z)\sigma=0$ by~\eqref{eq:Zsigma0}.
Thus the curvature of $\nabla^\Sigma|_{\Sigma_+}$ along $\mathcal{V}_+$ vanishes.
Since the Bott connection is also flat, we have $F_{\mathcal D^+}=F_{\nabla^\Sigma|_{\Sigma_+}}\otimes\Id+\Id\otimes F_{\nabla^B}=0$.
The flatness of $\mathcal{D}^-$ follows by Lemma \ref{lemma:conjugate}.
\end{proof}

Now we focus on transversal geometry.
The concepts and terminology from transversal geometry used in this paper are summarized in Appendix~\ref{app:transversal holomorphic geometry}.

\begin{proposition}\label{prop:basic holomorphic line bundle}
Let $(M^3,\phi,\xi,\eta,g)$ be a three-dimensional connected trans-Sasakian $\eta$-Einstein manifold.
Let $\mathcal F_\xi$ be the foliation by the flow lines of $\xi$.
Then $\mathcal F_\xi$ is a transversely holomorphic foliation whose transversal complex structure is induced by $\phi$.
Moreover, $\mathcal D^+$ defines a basic holomorphic line bundle structure on $\mathcal R_\xi^+$, and $\mathcal D^-$ defines a basic anti-holomorphic line bundle structure on $\mathcal R_\xi^-$.
\end{proposition}
\begin{proof}
Take flow boxes $U_i$ with projections $p_i:U_i\to S_i$ along the flow of $\xi$, here $S_i$ is a transversal surface in $U_i$.
By~\eqref{eq:L xi phi 0}, $\phi$ induces a complex structure on $S_i$.
The local maps between the surfaces $S_i$ obtained by following the flow lines preserve these complex structures and are therefore holomorphic, so $\mathcal F_\xi$ is a transversely holomorphic foliation.
Since $\mathcal D^+$ is flat by Proposition~\ref{prop:intrinsic-rs-flatness}, Rawnsley's theorem as recalled in~\cite[Section~2.2]{BiswasKasuya2023} gives a nowhere-vanishing $\mathcal D^+$-parallel frame $s_i$ of $\mathcal R_\xi^+$ on a sufficiently small flow box $U_i$.
On overlaps, write $s_j=g_{ij}s_i$.
For $X\in\Gamma(\mathcal{V}_+)$, we have $0=\mathcal D^+_X s_j=\mathcal D^+_X(g_{ij}s_i)=(Xg_{ij})s_i$, so $\xi g_{ij}=0$ and $\dd g_{ij}|_{\mathcal H^{0,1}}=0$.
Thus this system of parallel frames defines a basic holomorphic line bundle structure on $\mathcal R_\xi^+$.
The case of $\mathcal R_\xi^-$ is proved similarly.
\end{proof}

As shown in Proposition~\ref{prop:basic holomorphic line bundle}, $\mathcal{D}^\pm$-parallel sections of $\mathcal{R}_\xi^\pm$ exist at least locally, but their global existence is not trivial.
If such a section exists globally, it defines a global (anti-)holomorphic section of the basic (anti-)holomorphic line bundle and also a global RS field.
We denote the set of $\mathcal D^\pm$-parallel sections of $\mathcal R_\xi^\pm$ by \textcolor{blue}{\textbf{$H_B^0(\mathcal R_\xi^\pm)$}}.
Also, let \textcolor{blue}{\textbf{$\operatorname{RS}(M,g)$}} be the complex vector space of global RS fields on $M$, and set $\textcolor{blue}{\operatorname{RS}_\xi(M,g)}:=\{\Psi\in\operatorname{RS}(M,g):\Psi(\xi)=0\}$.
By Theorem~\ref{thm:fixed} and Proposition~\ref{prop:intrinsic-rs-connection}, we obtain the following theorem immediately.

\begin{theorem}\label{thm:global-rs-basic}
Let $(M^3,\phi,\xi,\eta,g)$ be a three-dimensional connected trans-Sasakian $\eta$-Einstein manifold.
Then there is a complex linear isomorphism $\operatorname{RS}_\xi(M,g)\simeq H_B^0(\mathcal R_\xi^+)\oplus H_B^0(\mathcal R_\xi^-)$.
Moreover, if $g$ does not have constant curvature, then $\operatorname{RS}(M,g)=\operatorname{RS}_\xi(M,g)$.
\end{theorem}

\section{Rank-two Rarita-Schwinger fields}
\label{sec:explicit-model-fields}
In this section, we give explicit expressions for all RS fields of type (R2) on three-dimensional simply connected and geodesically complete trans-Sasakian $\eta$-Einstein manifolds.

Let $(M,\phi,\xi,\eta,g)$ be a geodesically complete and simply connected three-dimensional trans-Sasakian $\eta$-Einstein manifold of type $(\alpha,\beta)$.
The ACM structure is $\alpha$-Sasakian when $\alpha\ne0,\beta\equiv0$, $\beta$-Kenmotsu when $\alpha\equiv0$, and proper trans-Sasakian otherwise, that is, when $\alpha\beta\not\equiv0$.
In~\cite{Matsuno2026}, we showed that every connected, simply connected, noncompact three-dimensional proper trans-Sasakian $\eta$-Einstein manifold with complete characteristic vector field is isomorphic to one of the model spaces $\mathcal{M}_K(v,s_0),\ \mathcal{M}_K(U)$.
We review these model spaces below.
The manifolds considered below are listed in the following proposition.

\begin{proposition}\label{prop:classification of 3dim trans sasakian eta-einstein}
Up to homothety, a geodesically complete and simply connected three-dimensional trans-Sasakian $\eta$-Einstein manifold is one of the following.
\end{proposition}
\begin{enumerate}
  \item Proper trans-Sasakian $\eta$-Einstein manifolds that are not Einstein.
  \item Sasakian manifolds diffeomorphic to $\R^3$.
  \item $\beta$-Kenmotsu $\eta$-Einstein manifolds diffeomorphic to $\R\times S^2$ or $\R\times \R^2$.
  \item Trans-Sasakian structures on $\mathbb{E}^3$.
  \item Two types of trans-Sasakian structures on $\mathbb{H}^3$.
  \item Sasakian manifolds diffeomorphic to $S^3$.
\end{enumerate}
\begin{proof}
If $g$ does not have constant curvature and $M$ is compact, then $M$ is homeomorphic to $S^3$, and~\cite{Matsuno2025} implies that the ACM structure is either $\alpha$-Sasakian or $\beta$-Kenmotsu.
Since $S^3$ admits no $\beta$-Kenmotsu structure, the ACM structure is $\alpha$-Sasakian.
Thus the $\eta$-Einstein condition and~\cite{Matsuno2025} imply that $M$ is homothetic to a Sasakian manifold.

Suppose that $g$ does not have constant curvature and that $M$ is noncompact.
If the ACM structure is proper trans-Sasakian, then $M$ is isomorphic as an almost contact metric manifold to one of the model spaces $\mathcal{M}_K(v,s_0),\ \mathcal{M}_K(U)$ by~\cite{Matsuno2026}.
If the ACM structure is $\alpha$-Sasakian, then $M$ is homothetic to a Sasakian manifold homeomorphic to $\R^3$ by~\cite{Matsuno2026}.
If the ACM structure is $\beta$-Kenmotsu, then $M$ is homeomorphic to $\R\times S^2$ or $\R\times\R^2$.

If $g$ has constant curvature, completeness and simple connectivity imply that, up to homothety, $(M,g)$ is $\mathbb E^3$, $\mathbb H^3$, or the round sphere $S^3$.
In each case, the $\eta$-Einstein condition is automatically satisfied.
The characteristic vector field generates a conformal foliation by geodesics~\cite[Proposition~3.3]{Mukud}, so the classification of Baird--Wood~\cite{BairdWood1991} applies.
On $S^3$, we obtain the Sasakian structure corresponding to the Hopf fibration.
On $\mathbb E^3$, we obtain the structure corresponding to a foliation by mutually parallel lines, while on $\mathbb H^3$, we obtain two inequivalent types of structures.
\end{proof}

\subsection{On proper trans-Sasakian manifolds}
We summarize below the definitions of the model spaces from~\cite{Matsuno2026}.
Let $K\in\R$ be a constant, and let $\D=\{z\in\C:|z|<1\}$ be the unit disk with its standard complex structure.
Let $z=x+\ii y$ and $\Delta_0=\partial_x^2+\partial_y^2$.
An \textbf{admissible pair} for $K>0$ is a pair $(v,s_0)\in C^\infty(\D,\R)\times\mathcal O(\D)$ satisfying
\begin{equation}
 \Delta_0v>0,\qquad e^v|s_0|<1,\qquad s_0\not\equiv0.
\end{equation}
An \textbf{admissible potential} for $K=0$ is a function $U\in C^\infty(\D,\R)$ satisfying
\begin{align}
 U>0,\qquad \Delta_0U<0.
\end{align}
An \textbf{admissible potential} for $K<0$ is a function $U\in C^\infty(\D,\R)$ satisfying
\begin{align}
 U\not\equiv0,\qquad |U|<\frac{\pi}{\sqrt{2|K|}},
\end{align}
and for which there exists a function $H_U\in C^\infty(\D,(0,\infty))$ satisfying
\begin{align}
 \Delta_0U=\sqrt{2|K|}\,H_U\sin(\sqrt{2|K|}U).
\end{align}
We denote the respective sets of data by $\mathcal U_K^{\mathrm{pair}}$ for $K>0$ and $\mathcal U_K^{\mathrm{pot}}$ for $K\leq0$.

For each fixed $K\in\R$, take an admissible pair $(v,s_0)$ if $K>0$, or an admissible potential $U$ if $K\leq0$, on the unit disk $\D$.
On $\D_z\times\R_w$, we define the one-form $\eta_K$ as
\begin{equation}
 \eta_K=
 \begin{cases}
  \dd w-(2K)^{-1/2}*_0\dd v,&K>0,\\
  \dd w+*_0\dd U,&K=0,\\
  \dd w-*_0\dd U,&K<0.
 \end{cases}
\end{equation}
On $\D_z\times\R_w$, we define the metrics as follows.
If $K>0$,
\begin{equation}
 \begin{aligned}
 &g_K=A_{K,v,s_0}|\dd z|^2+\eta_K^2,\\
 &A_{K,v,s_0}=\frac{\Delta_0v}{2K}\frac{|1-\zeta|^2}{1-|\zeta|^2},
 \qquad \zeta=e^{v-\ii\sqrt{2K}w}\overline{s_0}.
 \end{aligned}
\end{equation}
If $K=0$,
\begin{equation}
\begin{aligned}
&g_0=A_{0,U}|\dd z|^2+\eta_0^2,\\
&A_{0,U}=-\frac{\Delta_0U}{2U}(w^2+U^2),
\end{aligned}
\end{equation}
If $K<0$,
\begin{equation}
 \begin{aligned}
 &g_K=A_{K,U}|\dd z|^2+\eta_K^2,\\
 &A_{K,U}=H_U\{\cosh(\sqrt{2|K|}w)+\cos(\sqrt{2|K|}U)\}.
 \end{aligned}
\end{equation}
Set $\xi=\partial_w$ and let $\varphi_K$ be the horizontal lift of the standard complex structure on $\D$, with $\varphi_K\xi=0$.
We denote these ACM models by $\mathcal M_K(v,s_0)$ for $K>0$ and $\mathcal M_K(U)$ for $K\leq0$.
These model spaces are not automatically geodesically complete, so suitable pairs or potentials must be chosen to make the metrics geodesically complete.
In fact, in~\cite{Matsuno2026}, we showed that there exist families of such pairs and potentials with infinitely many degrees of freedom.
Since this paper considers only the geodesically complete case, we assume that the pair or potential is chosen so that the metric is geodesically complete.
However, the following proposition also holds for any pair or potential.

\begin{proposition}
\label{prop:explicit-model-rs}
Let $M=\mathcal M_K(v,s_0),\ \mathcal M_K(U)$ be a complete three-dimensional proper trans-Sasakian manifold that is not Einstein.
Denote the metric coefficients $A_{K,v,s_0}$, $A_{0,U}$, and $A_{K,U}$ by $A$, and write $ \eta_K=\dd w+B_x\dd x+B_y\dd y$.
Take the global orthonormal frame and its dual coframe to be
\begin{equation}\label{eq:explicit-model-frame}
 \begin{aligned}
 e_1&=\xi=\partial_w,&
 e_2&=A^{-1/2}(\partial_x-B_x\partial_w),&
 e_3&=A^{-1/2}(\partial_y-B_y\partial_w),\\
 \theta^1&=\eta_K,&
 \theta^2&=A^{1/2}\dd x,&
 \theta^3&=A^{1/2}\dd y
 \end{aligned}
\end{equation}
and let $\chi_\pm$ be orthonormal spinors with constant components relative to a spin lift of this frame, satisfying $\gamma_1\chi_\pm=\pm i\chi_\pm$.
Set $c:=\sqrt{2|K|}$ and define the following functions.
\begin{equation}\label{eq:explicit-model-amplitudes}
 \begin{array}{c|c}
  &F_\pm\\ \hline
 K>0&
 e^{-v/4\mp icw/4}\sqrt{1-\Re\zeta\pm i\Im\zeta}\\[3pt]
 K=0&\sqrt{U\mp iw}\\[3pt]
 K<0&
 \sqrt{\cosh\!\left(\frac c2(w\mp iU)\right)}
 \end{array}
\end{equation}
$F_\pm$ are smooth and nowhere-vanishing on the whole model and satisfy $F_-=\overline{F_+}$.

Set $\Phi=A^{-3/4} \left\{h_+(z)F_+\chi_++\overline{h_-(z)}F_-\chi_-\right\}$ for any two holomorphic functions $h_\pm\in\mathcal{O}(\D)$.
Then $\Psi=\Phi\theta^2+\gamma_1\Phi\theta^3$ is a global RS field satisfying $\Psi(\xi)=0$.
Conversely, every $\Psi\in \operatorname{RS}(M,g)=\operatorname{RS}_\xi(M,g)$ is uniquely expressed by this formula in terms of a pair of holomorphic functions $(h_+,h_-)$.
\end{proposition}
\begin{proof}
We construct a global $\mathcal{D}^+$-parallel section of $\mathcal{R}_\xi^+$.
Let $S:=\{w=0\}$.
This is a global surface transverse to $\xi$.
We construct a suitable initial section on $S$ and extend it in the $w$-direction by $\mathcal{D}_\xi^+$-parallel transport to compute a global parallel section.

Set $u:=\frac12\log A,\allowbreak\ B_z:=\frac{B_x-iB_y}{2},\allowbreak\ B_{\bar z}:=\frac{B_x+iB_y}{2},\allowbreak\ \partial_z:=\frac{\partial_x-i\partial_y}{2},\allowbreak\  \partial_{\bar z}:=\frac{\partial_x+i\partial_y}{2}$.
Take the complex frame $\xi=\partial_w,\ \partial=\sqrt2\,A^{-1/2}(\partial_z-B_z\partial_w),\ \bar\partial=\sqrt2\,A^{-1/2}(\partial_{\bar z}-B_{\bar z}\partial_w)$.
Then $\mathcal{V}_+=\operatorname{span}_\C\{\partial_w,\partial_{\bar z}\}$.
The functions $\alpha,\beta$ are given by $2\beta=\diver\xi=A^{-1}\partial_wA=2u_w$ and $\dd\eta=(\partial_xB_y-\partial_yB_x)\,\dd x\wedge\dd y=-2\alpha A\,\dd x\wedge\dd y$.
Set $\rho:=\beta+i\alpha$.

\paragraph{Computing the connection $\mathcal{D}^+$:}
We compute $\mathcal{D}^+$ using the global section $e:=A^{-3/4}\chi_+(\theta^2+i\theta^3)=e^{-u/2}\chi_+\dd z$ of $\mathcal{R}_\xi^+$.
The real orthonormal frame satisfies $[\xi,e_2]=-\beta e_2,\ [\xi,e_3]=-\beta e_3,\ [e_2,e_3]=2\alpha\xi+(e_3u)e_2-(e_2u)e_3$.
Set $\omega_{23}(X):=g(\nabla^g_Xe_2,e_3)$.
Koszul's formula gives $2\omega_{23}(\xi)=g([\xi,e_2],e_3)-g([e_2,e_3],\xi)+g([e_3,\xi],e_2)=-2\alpha,\ \omega_{23}(e_2)=-g([e_2,e_3],e_2)=-e_3u$, and $\omega_{23}(e_3)=-g([e_2,e_3],e_3)=e_2u$.
Since $\partial_{\bar z}=B_{\bar z}\xi+\frac{e^u}{2}(e_2+i e_3)$, we obtain $\omega_{23}(\partial_{\bar z})=-\alpha B_{\bar z}+\frac{e^u}{2}(-e_3u+i e_2u)=-\alpha B_{\bar z}+\frac12\bigl(-(u_y-B_yu_w)+i(u_x-B_xu_w)\bigr)=-\alpha B_{\bar z}+i(u_{\bar z}-B_{\bar z}u_w)$.
For $X\in\mathcal V_+$, we have $\nabla^g_X\xi\in\mathcal H^{0,1}$, so~\eqref{eq:Zsigma0} gives $\frac12g(\nabla^g_X\xi,e_2)\gamma_1\gamma_2\chi_++\frac12g(\nabla^g_X\xi,e_3)\gamma_1\gamma_3\chi_+=\frac12\gamma_1\gamma(\nabla^g_X\xi)\chi_+=0$.
Using $\gamma_2\gamma_3\chi_+=-i\chi_+$, we obtain $\nabla^\Sigma_X\chi_+=\frac12\gamma_1\gamma(\nabla^g_X\xi)\chi_++\frac12\omega_{23}(X)\gamma_2\gamma_3\chi_+=-\frac{i}{2}\omega_{23}(X)\chi_+$.
Finally, $\nabla^B_\xi\dd z=\nabla^B_{\partial_{\bar z}}\dd z=0$ and $\bar\rho=u_w-i\alpha$ give
\begin{align*}
\mathcal D^+_\xi e
&=\left(-\frac12u_w-\frac{i}{2}\omega_{23}(\xi)\right)e
=\left(-\frac12u_w+\frac{i\alpha}{2}\right)e
=-\frac12\bar\rho\,e,\\
\mathcal D^+_{\partial_{\bar z}}e
&=\left(-\frac12u_{\bar z}-\frac{i}{2}\omega_{23}(\partial_{\bar z})\right)e=\left(-\frac12u_{\bar z}+\frac{i\alpha}{2}B_{\bar z}+\frac12(u_{\bar z}-B_{\bar z}u_w)\right)e=-\frac12B_{\bar z}\bar\rho\,e.
\end{align*}

\paragraph{Computing the initial section on $S$ and the evolution equation:}
On $S$, set $\bar{s}:=e^q e$.
To satisfy the parallelism condition $\mathcal D^+_{\partial_{\bar z}}\bar{s}=0$, we must solve
\begin{align*}
\partial_{\bar z}q=\frac12B_{\bar z}\bar\rho|_{w=0}.
\end{align*}
We extend the initial section $\bar{s}=e^q e$ on $S$ obtained from a solution $q$ of this equation by parallel transport in the $\xi=\partial_w$ direction.
Writing the extension as $\Psi_+=F_+ e$, we obtain $\mathcal D_\xi(F_+ e) =(\partial_wF_+)e+F_+\mathcal D_\xi e =\left(\partial_wF_+-\frac12\bar\rho F_+\right)e$, so parallelism and the initial condition $\Psi_+|_S=\bar{s}$ give the equation $\partial_wF_+=\frac12\bar\rho F_+,\  F_+(z,0)=e^{q(z)}$.
Integration of this equation gives the following expression.
\begin{align*}
F_+(z,w):=\exp\!\left\{q(z)+\frac12\int_0^w\bar\rho(z,t)\,\dd t\right\}
\end{align*}
Moreover, since $\mathcal D^-=\overline{\mathcal D^+}$ and $F_-=\overline{F_+}$, a nowhere-vanishing $\mathcal{D}^-$-parallel section of $\mathcal{R}_\xi^-$ is given by $\Psi_-=A^{-3/4}F_-\chi_-(\theta^2-i\theta^3)$.
Thus we compute $q$ and $F_+$ below.

\paragraph{The case $K>0$:}
Since $B_{\bar z}=-\frac{i}{c}v_{\bar z},\ \bar\rho=-\frac{ic}{2}\frac{1+e^{v+icw}s_0}{1-e^{v+icw}s_0}$ and $s_0$ is holomorphic, we obtain
\begin{align*}
 \frac12B_{\bar z}\bar\rho|_{w=0}
 &=-\frac{v_{\bar z}}4\frac{1+e^v s_0}{1-e^v s_0}
 =\partial_{\bar z}
     \left\{-\frac v4+\frac12\Log(1-e^v s_0)\right\}.\\
\therefore\ q&=-\frac v4+\frac12\Log(1-e^v s_0).
\end{align*}
Next, since $\frac12\bar\rho=-\frac{ic}{4}      +\frac12\partial_w\Log(1-e^{v+icw}s_0)$, we obtain
\begin{align*}
 \frac12\int_0^w\bar\rho(z,t)\,\dd t
 &=-\frac{icw}{4}
       +\frac12\Log(1-e^{v+icw}s_0)
       -\frac12\Log(1-e^v s_0),\\
 \therefore\ F_+(z,w)
 &=\exp\!\left\{-\frac v4-\frac{icw}{4}
                  +\frac12\Log(1-e^{v+icw}s_0)\right\}
 =e^{-v/4-icw/4}\sqrt{1-\bar \zeta}.
\end{align*}

\paragraph{The case $K=0$:}
Since $B_{\bar z}=iU_{\bar z},\ \bar\rho=\frac1{w+iU}$, the same procedure gives the following expressions.
\begin{align*}
q=\frac12\log U,\ F_+(z,w)=\sqrt{U-iw}
\end{align*}

\paragraph{The case $K<0$:}
Noting that $\Delta_0U=cH_U\sin(cU)$, we have $B_{\bar z}=-iU_{\bar z},\ \bar\rho=\frac c2\frac{\sinh(cw)-i\sin(cU)}{\cosh(cw)+\cos(cU)}=\frac c2\tanh\!\left(\frac c2(w-iU)\right)$.
Thus a similar calculation gives the following expressions.
\begin{align*}
q=\frac12\log\cos\left(\frac{cU}{2}\right),\ F_+(z,w)=\sqrt{\cosh\!\left(\frac c2(w-iU)\right)}
\end{align*}
\end{proof}

\subsection{\texorpdfstring{On Sasakian manifolds diffeomorphic to $\R^3$}{On Sasakian manifolds diffeomorphic to R3}}
We give model spaces for Sasakian manifolds diffeomorphic to $\R^3$ as follows.
For $B=\C$ or $\D$, let $z=x+iy$ and let $v\in C^\infty(B,\R)$ satisfy $\Delta_0v>0$.
On $B_z\times\R_w\simeq\R^3$, we define a contact form $\eta$, a Reeb vector field $\xi$ and a metric $g$
\begin{align}
\eta=\dd w-\frac12*_0\dd v,\ 
\xi=\partial_w,\\
g=A|\dd z|^2+\eta^2,\ A=\frac14\Delta_0v. \label{eq:sasakian-r3-model}
\end{align}
Choose the potential $v$ so that the metric $g$ is geodesically complete.
We define the model space by $\mathcal{S}:=(\R^3,\phi,\xi,\eta,g)$, here $\phi(X):=\xi\times X$.
The model $\mathcal{S}$ is obtained from the proper trans-Sasakian $\eta$-Einstein space $\mathcal{M}_K(v,s_0)$ for $K>0$ by setting $s_0\equiv0,\ K=2$ and allowing $B=\D$ or $\C$.
Although $s_0\equiv0$ and $B=\C$ are excluded from the definition of $\mathcal{M}_K(v,s_0)$, the calculation in Proposition~\ref{prop:explicit-model-rs} remains valid with these formal choices.

\begin{proposition}\label{prop:sasakian-r3-model}
Every geodesically complete three-dimensional Sasakian manifold diffeomorphic to $\R^3$ is isomorphic to some model space $\mathcal{S}$.
Moreover, $\mathcal{S}$ does not have constant curvature.
\end{proposition}
\begin{proof}
Since the metric is geodesically complete, the unit Killing field $\xi$ is complete.
Since $M$ is connected, simply connected and noncompact, the complete nowhere-vanishing Killing field $\xi$ generates a free proper $\R$-action~\cite[Lemma~3.12]{Matsuno2026}.
Thus the orbit space $B=M/\R$ is a simply connected oriented surface, and we obtain a trivial principal bundle representation $M=B\times\R_t$.
Since the structure is Sasakian, we can further write $\xi=\partial_t,\ \eta=\dd t+\lambda,\ g=h+\eta^2,\ \dd\lambda=-2\operatorname{vol}_h$.
If $B$ were compact, applying Stokes' theorem to $\dd\lambda=-2\operatorname{vol}_h$ would give a contradiction, so $B$ is noncompact.
The tensor $\phi$ induces a complex structure on $B$, and the uniformization theorem allows us to take $B=\C$ or $\D$ and write $h=A|\dd z|^2$.
The Poisson equation $\Delta_0v=4A$ has a smooth real-valued solution on the whole of $B$~\cite[Corollary~8.3 and the subsequent remark]{McMullenRS}.
Thus $\dd\left(\lambda+\frac12*_0\dd v\right)=\left(-2A+\frac12\Delta_0v\right)\dd x\wedge\dd y=0$, we can write $\lambda+\frac12*_0\dd v=\dd f$.
Setting $w=t+f$ gives~\eqref{eq:sasakian-r3-model}.
Moreover, $\Ric(\xi,\xi)=2$, so if the space has constant curvature, this curvature is positive.
Geodesic completeness then implies that the space is the round $S^3$, which contradicts the assumption.
\end{proof}
The above immediately gives the following proposition.
\begin{proposition}\label{prop:sasakian-r3-rs}
Let $\mathcal{S}=(\R^3,\phi,\xi,\eta,g)$ be a geodesically complete Sasakian manifold.
Then $\operatorname{RS}(\R^3,g)=\operatorname{RS}_\xi(\R^3,g)$, and all RS fields are given by the formula in Proposition~\ref{prop:explicit-model-rs} with $s_0\equiv0,\ K=2$, $B=\C$ or $\D$, and $F_\pm=e^{-v/4\mp iw/2}$.
\end{proposition}

\subsection{\texorpdfstring{On $\beta$-Kenmotsu manifolds}{On beta-Kenmotsu manifolds}}
Let $(M^3,\varphi,\xi,\eta,g)$ be a simply connected, geodesically complete three-dimensional $\beta$-Kenmotsu $\eta$-Einstein manifold.
In~\cite{Matsuno2026}, we show that $(M,g)$ is given by
\[
M=\mathbb{R}_w\times B,\ 
g=dw^2+e^{2\Omega(w)}h,\ 
\Omega\in C^\infty(\R,\R)
\]
where $(B,h)$ is a connected, simply connected, complete Riemannian surface, and $B\simeq\R^2$ or $B\simeq S^2$.
If $B\simeq\R^2$, this family of spaces contains Euclidean space and hyperbolic space.
Since we consider these spaces separately, we exclude the case of constant curvature in this subsection.

\subsubsection{\texorpdfstring{$B\simeq\R^2$ case}{B isomorphic to R2 case}}

By the uniformization theorem, we identify $B$ biholomorphically with $\C$ or $\D$ and take a global complex coordinate $z=x+iy$.
Writing the base metric as $h=e^{2u(z,\bar z)}|dz|^2,\ u\in C^\infty(B,\R)$, we have $g=dw^2+e^{2(\Omega(w)+u(z,\bar z))}(dx^2+dy^2)$.
Take the orthonormal coframe to be $\theta^1=\eta=dw,\ \theta^2=e^{\Omega+u}dx,\ \theta^3=e^{\Omega+u}dy$.
Fix a global spin lift of its dual frame $\{e_1=\xi,e_2,e_3\}$.
Let $\chi_\pm$ be orthonormal spinors with constant components relative to this spin frame, satisfying $\gamma(\xi)\chi_\pm=\pm i\chi_\pm$.

\begin{proposition}\label{prop:kenmotsu-plane-rs}
Let $(\R^3,\varphi,\xi,\eta,g)$ be a complete three-dimensional $\beta$-Kenmotsu manifold given above which does not have constant curvature.
Then $\operatorname{RS}(\R^3,g)=\operatorname{RS}_\xi(\R^3,g)$, and every $\Psi\in \operatorname{RS}(\R^3,g)$ is given by
\begin{equation}\label{eq:kenmotsu-plane-rs}
 \Psi=e^{-\Omega-\frac32u}
 \left\{
 h_+(z)\chi_+(\theta^2+i\theta^3)
 +\overline{h_-(z)}\chi_-(\theta^2-i\theta^3)
 \right\},
\end{equation}
where $h_+,h_-\in\mathcal O(B)$ are uniquely determined holomorphic functions.
\end{proposition}
\begin{proof}
We proceed as in the proper trans-Sasakian case.
Since $\nabla_X\xi=\Omega'(w)\{X-\eta(X)\xi\}$, we have $\alpha=0,\ \beta=\Omega'(w)$.
We have $\mathcal{V}_+=\operatorname{span}_\C\{\partial_w,\partial_{\bar z}\}$.
Since $A=e^{2(\Omega+u)},\ B_x=B_y=0$, we obtain $ \nabla^\Sigma_{\partial_w}\chi_+=0,\ \nabla^\Sigma_{\partial_{\bar z}}\chi_+=\frac12u_{\bar z}\chi_+$.
Also, $\nabla^B_{\partial_w}dz=\Lie_{\partial_w}dz=0$ and $\nabla^B_{\partial_{\bar z}}dz=\Lie_{\partial_{\bar z}}dz=0$.
Taking the global section of $\mathcal{R}_\xi^+$ to be $e_+:=e^{-\frac32(\Omega+u)}  \chi_+(\theta^2+i\theta^3)=e^{-\frac12(\Omega+u)}\chi_+dz$, we find that $\mathcal D^+_{\partial_w}e_+=-\frac{\Omega'}2e_+,\ \mathcal D^+_{\partial_{\bar z}}e_+=e^{-\frac12(\Omega+u)}\left(-\frac12u_{\bar z}+\frac12u_{\bar z}\right)\chi_+dz=0$.
If we write the desired section as $\Psi_+=F_+e_+$, its parallelism is equivalent to $ \partial_wF_+=\frac{\Omega'}2F_+,\ \partial_{\bar z}F_+=0$, and solving these equations gives $F_+(z,w)=e^{\Omega(w)/2}h_+(z),\ h_+\in\mathcal O(B)$.
Since the $\mathcal{R}_\xi^-$ component is given by complex conjugation, we obtain~\eqref{eq:kenmotsu-plane-rs}.
\end{proof}

\subsubsection{\texorpdfstring{$B\simeq S^2$ case}{B isomorphic to S2 case}}
\begin{proposition}\label{prop:kenmotsu-sphere-vanishing}
In the case $B\simeq S^2$, no complete three-dimensional $\beta$-Kenmotsu $\eta$-Einstein manifold admits a nonzero RS field.
\end{proposition}
\begin{proof}
Since a complete simply connected three-dimensional  constant curvature space is one of $S^3$, $\R^3$, or $\mathbb H^3$, and none of these is diffeomorphic to $\R\times S^2$, the metric $g$ does not have constant curvature.
Thus $\operatorname{RS}(\R\times S^2,g)=\operatorname{RS}_\xi(\R\times S^2,g)$.
We identify the transverse surface $B_0:=\{w=0\}$ with $B=S^2$.

In each holomorphic coordinate $z$, write $h=e^{2u}|\dd z|^2$ and let the local frame from the previous case be $e_+=e^{-\frac32(\Omega+u)}\chi_+(\theta^2+i\theta^3)$.
We determine how $e_+$ transforms under a change of holomorphic coordinates.
For a holomorphic coordinate change $z'=f(z)$, write $f':=\dd f/\dd z$ and $f'/|f'|=e^{i\vartheta}$.
Then $u'=u-\log|f'|,\ \theta'^2+i\theta'^3=e^{\Omega+u'}\dd z'=\frac{f'}{|f'|}(\theta^2+i\theta^3),\ \chi'_+=\pm\exp\!\left(-\frac\vartheta2\gamma_2\gamma_3\right)\chi_+ =\pm\exp\!\left(\frac\vartheta2\gamma_1\right)\chi_+ =\pm e^{i\vartheta/2}\chi_+$, so the frame transforms as $e'_+=e^{-\frac32(\Omega+u')}\chi'_+(\theta'^2+i\theta'^3)=\pm|f'|^{3/2}e^{i\vartheta/2}e^{i\vartheta}e_+=\pm(f')^{3/2}e_+$.
Moreover, Proposition~\ref{prop:basic holomorphic line bundle} shows that $\mathcal{R}_\xi^+|_{B_0}$ is a holomorphic line bundle.

Let $K_B$ be the holomorphic cotangent bundle of $B=S^2$.
The above transformation law and $\dd z'=f'\dd z$ give $e'_+\otimes e'_+ =(f')^3(e_+\otimes e_+),\ \dd z'\otimes\dd z'\otimes\dd z'=(f')^3(\dd z\otimes\dd z\otimes\dd z)$, so the transition functions agree and we obtain a holomorphic bundle isomorphism $(\mathcal R_\xi^+|_{B_0})\otimes(\mathcal R_\xi^+|_{B_0})\simeq K_B\otimes K_B\otimes K_B$.

Since $B\simeq\mathbb{CP}^1$ and it is known that $K_B\simeq\mathcal O(-2)$, we have $ 2\deg(\mathcal R_\xi^+|_{B_0})=3\deg K_B=-6$.
Therefore, $\mathcal R_\xi^+|_{B_0}\simeq\mathcal O(-3)$.
If $\mathcal{R}_\xi^+|_{B_0}$ has a nonzero holomorphic section, the degree of the bundle equals the sum of the orders of the zeros of that section and must therefore be nonnegative.
Thus $\mathcal R_\xi^+|_{B_0}$, whose degree is negative, has no nonzero global holomorphic section.
Since $\mathcal D^+_{\partial_{\bar z}}e_+=0$, writing $\Psi_+|_{B_0}=a e_+$ locally gives $\partial_{\bar z}a=0$, so $a$ is holomorphic.
Hence $\Psi_+|_{B_0}$ is a global holomorphic section and therefore vanishes.
Since the $\mathcal{R}_\xi^-$ component is given by complex conjugation, its global anti-holomorphic sections also vanish.
\end{proof}

\subsection{\texorpdfstring{On Sasakian manifolds diffeomorphic to $S^3$}{On Sasakian manifolds diffeomorphic to S3}}\label{subsec:S3}
If $R^{T^*M}$ denotes the curvature of the cotangent bundle, the twisted Lichnerowicz formula~\cite[Chapter~II, Theorem~8.17]{LawsonMichelsohn1989} is $D_{TM}^2\Psi=\nabla^*\nabla\Psi+\frac{\Scal}{4}\Psi+\sum_{i<j}\gamma_i\gamma_jR^{T^*}_{e_i,e_j}\Psi$.
On the round $S^3$, if $\Pi\Psi=0$, then $D_{TM}^2\Psi=\nabla^*\nabla\Psi+\frac5{12}\Scal\,\Psi$.
Thus, for a global RS field, we obtain $0=\|\nabla\Psi\|_{L^2}^2+\frac5{12}\Scal\,\|\Psi\|_{L^2}^2$, so $\Psi=0$ because $\Scal>0$.
Therefore, in what follows, we consider the case in which $g$ does not have constant curvature.

Sasakian structures on $S^3$ have been classified~\cite[Sections~2.2 and~4]{WangZhang2013}.
Every Sasakian structure on $S^3$ is obtained by a second type deformation of a weighted Sasakian structure.
We give a brief review below.
\begin{definition}
Let $S^3:=\{(z_1,z_2)\in \C^2;\ |z_1|^2+|z_2|^2=1\}$.
For a pair $a=(a_1,a_2)$ of positive constants, the Reeb flow is $\varphi_t(z_1,z_2)=(e^{ia_1t}z_1,e^{ia_2t}z_2)$, and the corresponding Reeb vector field is $\xi_a=\sum_{j=1}^{2}ia_j(z_j\partial_{z_j}-\bar z_j\partial_{\bar z_j})$.
Set $\eta_0:=\sum_{i=1}^2\Im(\bar{z}_i\dd z_i),\ \sigma=\sum_j a_j|z_j|^2$, define the contact form by $\eta_a:=\eta_0/\sigma$, and let $\mathcal{H}_a=\ker\eta_a$.
Let $g_0$ be the standard round metric of sectional curvature $1$, and define the metric $g_a$ by $g_a(\xi_a,\xi_a)=1,\ g_a(\xi_a,\mathcal{H}_a)=0,\ g_a|_{\mathcal{H}_a}=\frac{1}{\sigma}g_0|_{\mathcal{H}_a}$.
Define the endomorphism field $\phi_a$ by $\dd\eta_a(X,Y)=2g_a(X,\phi_a Y)$.
We call $(S^3,\phi_a,\xi_a,\eta_a,g_a)$ the weighted Sasakian structure with weights $a=(a_1,a_2)$.
\end{definition}
The standard Sasakian structure corresponds to $a_1=a_2=1$.
For $a_1=a_2>0$, the Reeb flow is regular and the metric is a Berger metric.
If $a_1/a_2\in\Q$, then $\xi_a$ is quasi-regular, and if $a_1/a_2\notin\Q$, then $\xi_a$ is irregular.
Using coordinates on $S^3$ adapted to its foliation by tori, given by $z_1=\cos r\,e^{i\theta_1},\ z_2=\sin r\,e^{i\theta_2},\  r\in(0,\frac{\pi}{2}),\ \theta_i\in(0,2\pi)$, we have $\xi_a=a_1\partial_{\theta_1}+a_2\partial_{\theta_2}$.
In these coordinates, the metric is given by
\begin{align*}
g_a=\eta_a^2+\frac{dr^2}{\sigma}+\frac{\cos^2r\sin^2r}{\sigma}\left(\frac{a_2\,d\theta_1-a_1\,d\theta_2}{\sigma}\right)^2.
\end{align*}

Setting $\mathcal{H}_0:=\ker\eta_0$, we have $\mathcal{H}_a=\mathcal{H}_0$.
The endomorphism $\phi_0$ is then given by $\dd\eta_0(X,Y)=2g_0(X,\phi_0 Y)$ and satisfies $\phi_0|_{\mathcal H_0}=-J_0|_{\mathcal H_0}$ and $\phi_0\xi_0=0$, where $J_0$ is the standard complex structure on $\C^2$.
This is the convention of Boyer-Galicki.
In general, the kernel directions of $\phi_0$ and $\phi_a$ are different, so $\phi_0\ne\phi_a$, but on $\mathcal{H}_a=\mathcal{H}_0$ we have $\phi_a=\phi_0$, and hence the $\pm i$ eigenspaces of $\phi_0$ and $\phi_a$ coincide.

The following is a \textcolor{blue}{\textbf{second type deformation}}, which preserves the Reeb vector field and the transverse complex structure while changing the transverse Kähler metric.
For $(S^3,\phi_a,\xi_a,\eta_a,g_a)$, a real smooth function $u$ satisfying $\xi_a u=0$ is called a \textcolor{blue}{\textbf{basic potential}}.
Define $d_B^c=\frac i2(\bar\partial_B-\partial_B)$ with respect to the transverse complex structure, and set the new contact form to be $\eta=\eta_a+d_B^c u$.
We require $u$ to satisfy the positivity condition $-\frac12\dd(\eta_a+d_B^cu)(X,\phi_aX)>0\ (0\ne X\in\mathcal H_a)$.
Since we can check that $\eta(\xi_a)=1,\ \iota_{\xi_a}\dd\eta=0$, the Reeb vector field remains $\xi=\xi_a$, but $\mathcal{H}:=\ker\eta$ changes.
The correspondence between $\mathcal{H}_0$ and $\mathcal{H}$ is given, for $X\in\mathcal{H}_0$, by $X_u:=X-d_B^cu(X)\xi_a\in\mathcal{H}$.
Under this correspondence, $\phi_aX$ corresponds to $\phi X_u=\phi_aX-d_B^cu(\phi_aX)\xi_a$, so the new endomorphism is given by $\phi=\phi_a-\xi_a\otimes(d_B^c u\circ\phi_a)$.
It is straightforward to check that this satisfies $\phi^2=-\operatorname{Id}+\eta\otimes\xi$.
Finally, define $g=\eta\otimes\eta-\frac12\dd\eta(\,\cdot\,,\phi\,\cdot\,)$.
Since $(\phi_a,\xi_a,\eta_a,g_a)$ is Sasakian, we have $\mathcal{L}_\xi g=0$, so $(\phi,\xi,\eta,g)$ is also K-contact, and hence Sasakian.
Also, since $\mathcal H^{0,1}=\{Z-d_B^cu(Z)\xi;\ Z\in\mathcal H_a^{0,1}\}$, $\mathcal{V}_+=\C\xi\oplus\mathcal H^{0,1}$ remains the same subbundle under a second type deformation.

\begin{proposition}\label{prop:sasakian-s3-rs}
Let $S^3:=\{(z_1,z_2)\in \C^2;\ |z_1|^2+|z_2|^2=1\}$, and let $(S^3,\phi,\xi,\eta,g)$ be a Sasakian structure of nonconstant curvature obtained by a second type deformation with basic potential $u$ of a weighted Sasakian structure with weights $a=(a_1,a_2)$.
Then the space of global RS fields is given by
\begin{align}
&\operatorname{RS}(S^3,g)=\bigoplus_{(j,k)\in\mathcal I_a}\left(\C\,\bar z_1^{\,j}\bar z_2^{\,k}\varepsilon    \oplus\C\,z_1^jz_2^k I(\overline\varepsilon)\right)\\
&\mathcal I_a:=\{(j,k)\in\Z_{\ge0}^2:(2j+3)a_1+(2k+3)a_2=1\}.\nonumber
\end{align}
Here $\varepsilon$ is a suitable nowhere-vanishing global section of $\mathcal R_\xi^+$, and $I$ is the isomorphism given in Lemma~\ref{lemma:conjugate}.
Thus, a nonzero global RS field exists if and only if $\mathcal I_a\ne\varnothing$, and $\dim_\C\operatorname{RS}(S^3,g)=2\#\mathcal I_a$.
\end{proposition}
\begin{proof}
\leavevmode\par
\begin{enumerate}
\item There exist a global complex vector field $\partial$ satisfying $\mathcal{V}_+=\operatorname{span}_{\C}\{\xi,\partial\}$ and a nowhere-vanishing global section $\lambda$ of $\operatorname{Ann}(\mathcal{V}_+)$.\\
\hspace*{2em} Proof : Since $\phi_a|_{\mathcal H_0}=\phi_0|_{\mathcal H_0}=-J_0|_{\mathcal H_0}$ for the standard complex structure $J_0$ on $\C^2$, we have $\mathcal H_a^{0,1}=\left.T^{1,0}\C^2\right|_{S^3}\cap(TS^3\otimes\C)$.
Hence we can define $\partial=-\bar z_2\partial_{z_1}+\bar z_1\partial_{z_2}$.
The annihilator of $\mathcal{V}_+=\C\xi\oplus\mathcal H^{0,1}=\C\xi\oplus\mathcal H^{0,1}_a$ is $\operatorname{Ann}(\mathcal{V}_+)=(\mathcal{H}^{1,0}_a)^*$, so we can take $\lambda=a_1\bar z_1\dd\bar z_2-a_2\bar z_2\dd\bar z_1$ as a global section.
In the weighted structure, $\partial$ spans $\mathcal H_a^{0,1}$, so $\mathcal{V}_+=\operatorname{span}_{\C}\{\xi,\partial\}$ also holds after a second type deformation.
Since $\lambda(\xi)=\lambda(\partial)=0$ and $\lambda(\bar\partial)=\sigma>0$, $\lambda$ is a nowhere-vanishing global frame of $\operatorname{Ann}(\mathcal{V}_+)$.

\item Set $\nu:=\frac{1-3(a_1+a_2)}2$.
There exists a nowhere-vanishing global section $\varepsilon$ of $\mathcal R_\xi^+$ satisfying $\mathcal D_{\partial}^+\varepsilon=0,\ \mathcal D_\xi^+\varepsilon=i\nu\varepsilon$.
For any $f\in C^\infty(S^3,\C)$, $f\varepsilon$ is a basic holomorphic section if and only if $\partial f=0,\ \xi f=-i\nu f$.\\
\hspace*{2em} Proof : Set $\mu=||\lambda||_g/\sqrt2$ and $\omega=\mu^{-1}\lambda=\theta^2+i\theta^3$.
Then $(\eta,\theta^2,\theta^3)$ is a global orthonormal coframe.
We denote by $(\xi,e_2,e_3)$ its dual frame.
Since $S^3$ is simply connected, this frame admits a global spin lift.
Choose a unit spinor $\chi_+$ with constant components relative to this lift satisfying $\gamma(\xi)\chi_+=i\chi_+$, and set $v=\mu^{3/2}\chi_+\otimes\omega$.
We compute $\mathcal{D}^+v$.
Set $\kappa(X)=g(\nabla_Xe_2,e_3)$.
For $X\in\mathcal{V}_+$, we have $\nabla_X\xi=-\phi X\in\mathcal H^{0,1}$, so equation~\eqref{eq:Zsigma0} gives $\gamma(\nabla_X\xi)\chi_+=0$.
Thus,
\begin{align*}
\nabla_X^\Sigma\chi_+&=\frac12\Bigl[g(\nabla_X\xi,e_2)\gamma_1\gamma_2+g(\nabla_X\xi,e_3)\gamma_1\gamma_3+\kappa(X)\gamma_2\gamma_3\Bigr]\chi_+\\
&=\frac12\gamma_1\gamma(\nabla_X\xi)\chi_++\frac12\kappa(X)\gamma_2\gamma_3\chi_+=-\frac i2\kappa(X)\chi_+.
\end{align*}
Since $\nabla_Y\xi=-\phi Y=\theta^3(Y)e_2-\theta^2(Y)e_3$, we have $\nabla_Ye_2=-g(e_2,\nabla_Y\xi)\xi+\kappa(Y)e_3=-\theta^3(Y)\xi+\kappa(Y)e_3$, and similarly $\nabla_Ye_3=\theta^2(Y)\xi-\kappa(Y)e_2$.
Hence $\nabla_Y\omega=g(\nabla_Ye_2,\cdot)+ig(\nabla_Ye_3,\cdot)=\{-\theta^3(Y)+i\theta^2(Y)\}\eta+\kappa(Y)(\theta^3-i\theta^2)=i\omega(Y)\eta-i\kappa(Y)\omega$.
Therefore, $d\omega(Y,Z)=(\nabla_Y\omega)(Z)-(\nabla_Z\omega)(Y)=i\{\omega(Y)\eta(Z)-\omega(Z)\eta(Y)\}-i\{\kappa(Y)\omega(Z)-\kappa(Z)\omega(Y)\}=-i\bigl((\kappa+\eta)\wedge\omega\bigr)(Y,Z)$.
For $X\in\mathcal{V}_+$, we have $\nabla_X^B\omega=\Lie_X\omega=\iota_X\dd\omega=-i(\kappa+\eta)(X)\omega$.
Since \(\mathcal{V}_+\) is involutive, for $X,Y\in\Gamma(\mathcal{V}_+)$ we have $(\Lie_X\lambda)(Y)=X\bigl(\lambda(Y)\bigr)-\lambda([X,Y])=0$.
Thus there is a unique smooth complex-valued function \(b(X)\) such that $\Lie_X\lambda=b(X)\lambda$.
Taking the Lie derivative of $\lambda=\mu\omega$ gives
\begin{align*}
 b(X)\lambda&=(X\mu)\omega+\mu\Lie_X\omega=\{X\log\mu-i(\kappa+\eta)(X)\}\lambda\\
 \therefore\ b(X)&=X\log\mu-i(\kappa+\eta)(X).
\end{align*}
Thus, for $X\in\mathcal{V}_+$,
\begin{align*}
\mathcal D_X^+v
 =\left\{\frac32X\log\mu-\frac{3i}{2}\kappa(X)-i\eta(X)\right\}v
 =\left\{\frac32b(X)+\frac i2\eta(X)\right\}v.
\end{align*}
We compute $b(\xi),\ b(\partial)$.
Since $\xi=ia_1\left(z_1\partial_{z_1}-\bar z_1\partial_{\bar z_1}\right) +ia_2\left(z_2\partial_{z_2}-\bar z_2\partial_{\bar z_2}\right)$ and $d\lambda=(a_1+a_2)\,d\bar z_1\wedge d\bar z_2$, we obtain $\Lie_\xi\lambda=d(\iota_\xi\lambda)+\iota_\xi d\lambda=-i(a_1+a_2)\lambda$ and $\Lie_{\partial}\lambda=d(\iota_{\partial}\lambda)+\iota_{\partial}d\lambda=0$.
Hence $b(\xi)=-i(a_1+a_2),\ b(\partial)=0$.
We compute $\eta(\partial)=\eta_a(\partial)+d_B^cu(\partial)=d_B^cu(\partial)$.
In a local transverse holomorphic coordinate $\zeta$, we have $d_B^cu=\frac{i}{2}(u_{\bar\zeta}\dd\bar{\zeta}-u_\zeta\dd\zeta)$ and $\partial=A\xi+B\partial_{\bar\zeta}$, so $\eta(\partial)=d_B^cu(\partial)=\frac{i}{2}\partial u$.
Thus,
\begin{align*}
\mathcal D_\xi^+v=i\nu v,\quad \mathcal D_{\partial}^+v=-\frac14(\partial u)v.
\end{align*}
Therefore,
\begin{align*}
\varepsilon:=e^{u/4}v=e^{u/4}\mu^{3/2}\chi_+\otimes\omega
\end{align*}
is a nowhere-vanishing global section satisfying $\mathcal D_{\partial}^+\varepsilon=0,\ \mathcal D_\xi^+\varepsilon=i\nu\varepsilon$.
Writing an arbitrary global section as $f\varepsilon$, we see that $f\varepsilon\in H_B^0(\mathcal R_\xi^+)$ is equivalent to $\partial f=0,\ \xi f=-i\nu f$.

\item $H_B^0(\mathcal R_\xi^+)$ has $\bar z_1^{\,j}\bar z_2^{\,k}\varepsilon,\ (j,k)\in\mathcal I_a$ as a basis over $\C$.\\
\hspace*{2em} Proof : In the coordinates $\{r,\theta_1,\theta_2\}$ given by $z_1=\cos r\,e^{i\theta_1},\ z_2=\sin r\,e^{i\theta_2},\ 0<r<\pi/2$, take the Fourier expansion
\[
 f=\sum_{j,k\in\Z}F_{jk}(r)e^{-i(j\theta_1+k\theta_2)}.
\]
Since $\xi=a_1\partial_{\theta_1}+a_2\partial_{\theta_2}$, the equation $\xi f=-i\nu f$ gives $(ja_1+ka_2-\nu)F_{jk}=0$.
On the other hand, since
\[
 \partial=\frac{e^{-i(\theta_1+\theta_2)}}2\left(\partial_r+i\tan r\,\partial_{\theta_1}-i\cot r\,\partial_{\theta_2}\right),
\]
the equation $\partial f=0$ gives $F_{jk}'+(j\tan r-k\cot r)F_{jk}=0$.
Thus,
\begin{align*}
F_{jk}=c_{jk}(\cos r)^j(\sin r)^k.
\end{align*}
Since $f$ is smooth on $S^3$, we have $|F_{jk}(r)|\le\|f\|_\infty$.
Hence, if $c_{jk}\ne0$, boundedness as $r\to\pi/2$ and $r\to0$ gives $j\ge0$ and $k\ge0$, respectively.
Moreover, $ja_1+ka_2=\nu$ is equivalent to $(j,k)\in\mathcal I_a$.
The Fourier expansion therefore becomes $f=\sum_{(j,k)\in\mathcal I_a}c_{jk}\bar z_1^{\,j}\bar z_2^{\,k}$.
Conversely, each monomial is smooth on all of $S^3$, and the identities $\partial\bar z_\ell=0$ and $\xi\bar z_\ell=-ia_\ell\bar z_\ell$ show that it satisfies both equations.
Since distinct monomials have distinct Fourier modes, they are linearly independent, and we obtain $H_B^0(\mathcal R_\xi^+)=\bigoplus_{(j,k)\in\mathcal I_a}\C\,\bar z_1^{\,j}\bar z_2^{\,k}\varepsilon$.
\end{enumerate}
\end{proof}

We remark that the above result depends only on the weights $a=(a_1,a_2)$ and not on the basic potential $u$ used in the second type deformation.
When $a_1=a_2=a$, the Reeb flow is regular, and the common minimal period of the fibers is $T=\frac{2\pi}{a}$.
In this case, the necessary and sufficient condition for the existence of an RS field can be written as $T=4\pi(d+3),\ d\in\mathbb Z_{\ge0}$.
In the quasi-regular case, that is, when $a_1=cp,\ a_2=cq,\ \textrm{gcd}(p,q)=1$, set $N:=\frac1{2c}-\frac32(p+q)$.
Then the necessary and sufficient condition for the existence of an RS field is $N\in p\mathbb Z_{\ge0}+q\mathbb Z_{\ge0}$.
In the irregular case, that is, when $a_1/a_2\not\in\mathbb{Q}$, either $\mathcal I_a=\varnothing$ or $\mathcal I_a=\{(j_0,k_0)\}$ for a unique pair of integers, and $\dim_{\mathbb C}\operatorname{RS}(S^3,g)\in\{0,2\}$.

In~\cite{Ohno2026LieGroups}, Ohno classified left-invariant RS fields on three-dimensional Lie groups and obtained left-invariant RS fields on a Berger sphere.
These fields are recovered from the above proposition in the case $u=0,\ a_1=a_2=1/6$.
In this case, $\mathcal I_a=\{(0,0)\}$, so $\operatorname{RS}(S^3,g)=\C\varepsilon\oplus\C I(\overline\varepsilon)$ and $\dim_\C\operatorname{RS}(S^3,g)=2$.

B\"ar and Mazzeo showed that, in dimension \(n\ge4\), for any given positive number, there exists a compact manifold on which the space of all RS fields has dimension greater than that number~\cite[Corollary~4.2]{BarMazzeo2021}.
As a corollary of Proposition~\ref{prop:sasakian-s3-rs}, we conclude that, for any given positive number, there exists a metric on $S^3$ for which the space of RS fields has dimension greater than that number.
In the construction of~\cite[Corollary~4.2]{BarMazzeo2021}, the dimension $n$ is fixed but the underlying topological manifold is not, whereas in the following corollary the underlying topological manifold can also be fixed to be $S^3$.
\begin{corollary}\label{coro:S3 regular dimRS}
Let $(S^3,\phi,\xi,\eta,g)$ be a Sasakian structure obtained by a second type deformation of the regular weighted Sasakian structure on $S^3$ with weights $a_1=a_2=\frac{1}{2(d+3)},\ d\in\mathbb Z_{\ge0}$.
Then $\dim_{\mathbb C}\operatorname{RS}(S^3,g)=2(d+1)$.

\end{corollary}
\begin{proof}
Since $\mathcal I_a=\{(j,k)\in\Z_{\ge0}^2:j+k=d\}=\{(j,d-j):j=0,\ldots,d\}$ in Proposition~\ref{prop:sasakian-s3-rs}, we obtain $\dim_\C\operatorname{RS}(S^3,g)=2\#\mathcal I_a=2(d+1)$.
\end{proof}

\subsection{\texorpdfstring{On $\mathbb{E}^3$ and $\mathbb{H}^3$}{On E3 and H3}}\label{subsec:E3 and H3}
In this subsection, we determine RS fields of type (R2) on $\mathbb E^3$ and $\mathbb H^3$.
Since a trans-Sasakian structure defines a conformal foliation by geodesics~\cite[Proposition~3.3]{Mukud}, we can use the classification of such foliations.
The classification of Baird--Wood~\cite[Theorems~2.14 and~3.3]{BairdWood1991} shows that every conformal foliation by geodesics defined on the whole of $\mathbb E^3$ consists of parallel lines.
On the whole of $\mathbb H^3$, there are two types, foliations by geodesics sharing a single ideal boundary point and foliations by geodesics orthogonal to a fixed totally geodesic plane $\mathbb H^2$.
These foliations can all be expressed as follows.
\begin{equation}\label{eq:spaceform-warped-model}
 M=\R_w\times B,\qquad
 g=\dd w^2+a(w)^2\lambda(z)^2|\dd z|^2,\qquad
 z=x+iy,\qquad \xi=\partial_w
\end{equation}
Here $B$, $a$, and $\lambda$ are given in the following table.
\begin{center}
\begin{tabular}{|l|p{0.5\linewidth}|l|l|l|}
\hline
Space & Foliation & $B$ & $a(w)$ & $\lambda(z)$\\
\hline
$\mathbb E^3$ & parallel lines & $\C$ & $1$ & $1$\\
\hline
$\mathbb H^3$ & geodesics with a common ideal boundary point & $\C$ & $e^w$ & $1$\\
\hline
\rule{0pt}{3.5ex}$\mathbb H^3$ & geodesics normal to a totally geodesic plane & $\D$ & $\cosh w$ & $\displaystyle\frac{2}{1-|z|^2}$\\
\hline
\end{tabular}
\end{center}
Take the orthonormal coframe to be $\theta^1=\dd w,\ \theta^2=a\lambda\,\dd x,\ \theta^3=a\lambda\,\dd y$.
Fix a global spin lift of its dual frame $\{e_1=\xi,e_2,e_3\}$.
Let $\chi_\pm$ be orthonormal spinors with constant components relative to this spin frame, satisfying $\gamma(\xi)\chi_\pm=\pm i\chi_\pm$.

\begin{proposition}\label{prop:complete-spaceform-rs}
Every RS field $\Psi$ of type (R2) on $\mathbb E^3$ or $\mathbb H^3$ can be written, after fixing a choice of $\xi$ satisfying $\Psi(\xi)=0$ and the above three coordinates adapted to $\xi$, as
\begin{equation}
 \Psi=\lambda(z)^{-1/2}
 \left\{h_+(z)\chi_+\,\dd z
       +\overline{h_-(z)}\chi_-\,\dd\bar z\right\},
 \qquad h_\pm\in\mathcal O(B).
\end{equation}
\end{proposition}
\begin{proof}
By Theorem~\ref{thm:removable-main} and the classification of Baird--Wood, the kernel foliation has one of the three standard forms above, and $\Psi(\xi)=0$ holds everywhere.
The calculation in Proposition~\ref{prop:kenmotsu-plane-rs} is also valid when $\Omega'\equiv0$, so we can perform the same calculation and obtain $F_+e_+ =e^{\Omega/2}h_+\,e^{-\frac32(\Omega+u)}\chi_+(\theta^2+i\theta^3)=\lambda^{-1/2}h_+\chi_+\,\dd z$.
Similarly, $F_-e_-=\lambda^{-1/2}\overline{h_-}\chi_-\,\dd\bar z$.
\end{proof}
Since these are spaces of constant curvature, not every RS field can be written in this form, and RS fields of type (R3) do in fact exist.

\section{Conclusions and discussions}\label{sec:conclusion}
In this paper, we showed that RS fields on three-dimensional Riemannian spin manifolds are of one of two types, namely, (R2) rank two almost everywhere and (R3) rank three almost everywhere, that the existence of an RS field of type (R2) induces a trans-Sasakian $\eta$-Einstein structure, and, in particular, that every RS field on trans-Sasakian $\eta$-Einstein manifold is of type (R2) if the metric does not have constant curvature.
This also shows that RS fields of types (R2) and (R3) cannot coexist on a manifold of nonconstant curvature.
We further showed that the RS equation reduces to a transversal holomorphic equation with respect to the characteristic foliation.
Using this reduction, we also classified RS fields of type (R2) on simply connected, geodesically complete three-dimensional Riemannian spin manifolds.
In the noncompact case, except for $M=S^2\times\R$, RS fields of type (R2) depend on arbitrary holomorphic functions on $\C$ or $\D$, so we have $\dim_{\C}\operatorname{RS}(M,g)=\infty$.
In the compact case, the round $S^3$ admits no nontrivial RS fields of either type (R2) or type (R3), but RS fields of type (R2) exist when the Sasakian structure satisfies a suitable condition analogous to Bohr's quantization condition.
Since $\ker D_{TM}$ on a compact manifold is finite-dimensional, the space of all RS fields is also finite-dimensional, but we showed that $\sup_{g}\dim_{\C}\operatorname{RS}(S^3,g)=\infty$.
This is the three-dimensional counterpart of the result of B\"ar and Mazzeo in dimensions four and higher.

Since we did not study RS fields of type (R3) in this paper, their classification is left for future work.
However, we obtain the following result on the dimension of the space of all RS fields of type (R3).

\begin{proposition}\label{prop:distinct-ricci-eigenvalues}
Let $(M,g)$ be a connected three-dimensional Riemannian spin manifold which does not have constant curvature.
If $\Psi$ is an RS field of type (R3), then we have $ \operatorname{RS}(M,g)=\C\Psi\oplus\C\,\mathfrak j\Psi$ and $\dim_{\C}\operatorname{RS}(M,g)=2$, where $\mathfrak j$ is the parallel quaternionic structure on $\Sigma M$.
\end{proposition}
\begin{proof}
See Appendix \ref{app:r3}.
\end{proof}

Compared with the infinite-dimensional space of RS fields of type (R2) on noncompact manifolds, the set of RS fields of type (R3) on manifolds of nonconstant curvature is small.
The space of all RS fields of type (R3) on the constant-curvature spaces $\mathbb{E}^3$ and $\mathbb{H}^3$ should also be infinite-dimensional.

\section*{Acknowledgements}
I am grateful to S.Ohno for discussions on the significance and novelty of this paper.
I also thank F.Ueno for helpful comments on the presentation of this paper.
During the preparation of this manuscript, AI was used to check proofs, search the literature for prior work, identify existing theorems applicable to the proofs, and proofread the English.
I conjectured the statements of the two lemmas proved in Appendices, and AI assisted me in developing their proofs.

This work was partly supported by MEXT Promotion of Distinctive Joint Research Center Program JPMXP0723833165 and Osaka Metropolitan University Strategic Research Promotion Project (Development of International Research Hubs).

\appendix
\section{Proof of Lemma \ref{lem:N elliptic-inequality}}\label{app:N elliptic-inequality}
Here we show that the following holds in the setting of the lemma.
\begin{equation}\label{eq:main}
\left|\nabla^*\nabla N+\frac54\Scal\,N\right|
\le C\left(
\frac{|\nabla\Psi|}{|\Psi|}|\nabla N|
+\frac{|\nabla\Psi|^2}{|\Psi|^2}|N|
\right)
\qquad\text{on }M^\circ.
\end{equation}
Here $C$ is a constant.
In particular, $|\nabla^*\nabla N|\le C_U(|N|+|\nabla N|)$ holds on every $U\Subset M^\circ$.

Let $L=\nabla^*\nabla,\ s=\Scal,\ Q=\Ric^\sharp$.
For a positively oriented local orthonormal frame $\{e_i\}$, set $\psi=(\psi_1,\psi_2,\psi_3)$ with $\psi_i=\Psi(e_i)$.
Set $T_{ai}:=(\nabla_{e_a}\Psi)(e_i),\ H(\psi,T):=\sum_a\bigl(T_{a1}\wedge T_{a2}\wedge\psi_3+T_{a1}\wedge\psi_2\wedge T_{a3}+\psi_1\wedge T_{a2}\wedge T_{a3}\bigr)$.

\paragraph{Step 1}
We have $LN=-\frac54sN-2*_\Sigma H(\psi,T)$.
\begin{proof}
With respect to an orthonormal frame $\{e_i\}$, the twisted Lichnerowicz formula~\cite[Chapter~II, Theorem~8.17]{LawsonMichelsohn1989} gives $(L\Psi)(X)=-\frac s4\Psi(X)+\sum_{j<k}\gamma_j\gamma_k\Psi\bigl(R(e_j,e_k)X\bigr)$.
For a frame satisfying $Qe_i=\rho_i e_i$ at a point, the curvature formula in three dimensions gives $R(e_i,e_j)e_i=-\left(\rho_i+\rho_j-\frac s2\right)e_j$.
We obtain from $0=\gamma_i\Pi(\Psi)=\sum_j\gamma_i\gamma_j\psi_j=-\psi_i+\sum_{j\ne i}\gamma_i\gamma_j\psi_j$ and $0=\gamma_i\Pi_{\mathrm{Ric}}(\Psi)=\sum_j\rho_j\gamma_i\gamma_j\psi_j=-\rho_i\psi_i+\sum_{j\ne i}\rho_j\gamma_i\gamma_j\psi_j$ that $(L\Psi)(e_i)
=-\frac s4\psi_i-\sum_{j\ne i}\left(\rho_i+\rho_j-\frac s2\right)\gamma_i\gamma_j\psi_j=\left(\frac s4-2\rho_i\right)\psi_i$, thus $L\Psi=\frac s4\Psi-2\Psi\circ Q$ holds.

Since $\nabla_{e_a}N=*_\Sigma\bigl(T_{a1}\wedge\psi_2\wedge\psi_3+\psi_1\wedge T_{a2}\wedge\psi_3+\psi_1\wedge\psi_2\wedge T_{a3}\bigr)$, we obtain $LN=*_\Sigma\Bigl[(L\Psi)(e_1)\wedge\psi_2\wedge\psi_3+\psi_1\wedge(L\Psi)(e_2)\wedge\psi_3+\psi_1\wedge\psi_2\wedge(L\Psi)(e_3)-2H(\psi,T)\Bigr] =\left\{\sum_{i=1}^3\left(\frac s4-2\rho_i\right)\right\}N-2*_\Sigma H(\psi,T)=-\frac54sN-2*_\Sigma H(\psi,T)$.
\end{proof}

Fix an arbitrary point $q\in M^\circ$ and set $\psi:=\Psi_q$.

\paragraph{Step 2}
There exist an orthonormal basis $\{e_i\}$ of $T_qM$, a nonzero spinor $\Phi\in\Sigma_qM$, and $t\in[0,1]$ such that $\psi=(t\gamma_3\Phi,\Phi,(1+t)\gamma_1\Phi)$ and the following estimates hold:
\begin{equation}\label{eq:scale}
\sqrt2|\Phi|\le|\psi|\le\sqrt6|\Phi|,
\qquad|\psi_1|\le |\Phi|,
\qquad|\psi_1|\le\frac{|N|}{|\Phi|^2}.
\end{equation}
\begin{proof}
The real Gram operator \(\psi^*\psi:T_qM\to T_qM\) is a self-adjoint operator satisfying $g\bigl((\psi^*\psi)X,Y\bigr)=\operatorname{Re}\langle\psi(X),\psi(Y)\rangle_\Sigma$.
We can choose an orthonormal eigenbasis for this operator such that $(\psi^*\psi)e_i=\lambda_i e_i,\ 0\le\lambda_1\le\lambda_2\le\lambda_3$.
If we set \(\psi_i=\psi(e_i)\), then $\operatorname{Re}\langle\psi_i,\psi_j\rangle_\Sigma=g\bigl((\psi^*\psi)e_i,e_j\bigr)=\lambda_i\delta_{ij}$ gives $|\psi_i|^2=\lambda_i$, and hence $\psi_i\perp\psi_j\ (i\ne j),\ |\psi_1|\le|\psi_2|\le|\psi_3|$.

Set $\Phi=\psi_2$.
If $\Phi=0$, then $\psi_1=0$, so $\Pi(\psi)=0$ gives $\psi_3=0$, and hence $\psi=0$.
Therefore, $\Phi\ne0$.
Since $\Phi,\gamma_1\Phi,\gamma_2\Phi,\gamma_3\Phi$ form a real orthogonal basis whose vectors have the same length, we can write $\psi_1=a\gamma_1\Phi+b\gamma_2\Phi+t\gamma_3\Phi$.
By $\Pi(\psi)=0$, we have $\psi_3=-\gamma_2\psi_1+\gamma_1\Phi=b\Phi+(1+t)\gamma_1\Phi-a\gamma_3\Phi$.
From $\psi_1\perp\psi_3$ and $\Phi\perp\psi_3$, we obtain $0=\Re\langle\psi_1,\psi_3\rangle=a|\Phi|^2,\ 0=\Re\langle\Phi,\psi_3\rangle=b|\Phi|^2$.
The ordering of the lengths gives $|t|\le1$ and $|1+t|\ge1$, so $0\le t\le1$.
Therefore,
\begin{equation}\label{eq:normalform}
\psi=(t\gamma_3\Phi,\ \Phi,\ (1+t)\gamma_1\Phi),\qquad0\le t\le1.
\end{equation}
Thus, \eqref{eq:normalform} gives $|\psi|^2=2(1+t+t^2)|\Phi|^2,\ |N|=t(1+t)|\Phi|^3$.
In particular, we have \eqref{eq:scale}.
\end{proof}

All calculations below are at the point $q$ fixed above, using the normal form~\eqref{eq:normalform}.
Set $E=\operatorname{span}_{\mathbb R}\{\Phi,\gamma_1\Phi\},\ P=T_{11},\ A=T_{12}$, and $B=T_{22}$.
Write $W=W_E+W_\perp$ for the real orthogonal decomposition of a spinor with respect to $E\oplus E^\perp$.

\paragraph{Step 3}
There is a constant $C'$ such that
\begin{equation}\label{eq:H-component-bound}
|H(\psi,T)|\le C'\bigl(|\psi_1||T|^2+|\Phi||T|(|P|+|A_\perp|)\bigr).
\end{equation}
\begin{proof}
Set $\psi_0=(0,\Phi,\gamma_1\Phi)$ for $t=0$.
Set $k=\psi_3-\gamma_1\Phi=-\gamma_2\psi_1$, so that $|k|=|\psi_1|$.
By the definition of $H$, we have $H(\psi,T)-H(\psi_0,T)=\sum_a\Bigl[T_{a1}\wedge T_{a2}\wedge\psi_3+T_{a1}\wedge\Phi\wedge T_{a3}+\psi_1\wedge T_{a2}\wedge T_{a3}-T_{a1}\wedge T_{a2}\wedge\gamma_1\Phi-T_{a1}\wedge\Phi\wedge T_{a3}\Bigr]=\sum_a\Bigl[T_{a1}\wedge T_{a2}\wedge(\psi_3-\gamma_1\Phi)+\psi_1\wedge T_{a2}\wedge T_{a3}\Bigr]=\sum_a\bigl[T_{a1}\wedge T_{a2}\wedge k+\psi_1\wedge T_{a2}\wedge T_{a3}\bigr]$.
It follows that $|H(\psi,T)-H(\psi_0,T)|\le|\psi_1|\sum_a\bigl(|T_{a1}||T_{a2}|+|T_{a2}||T_{a3}|\bigr)\le|\psi_1|\sum_a\left(\frac12|T_{a1}|^2+|T_{a2}|^2+\frac12|T_{a3}|^2\right)$.
Thus, we obtain
\begin{equation}\label{eq:H-difference1}
|H(\psi,T)-H(\psi_0,T)|\le|\psi_1||T|^2.
\end{equation}

Next we estimate $H(\psi_0,T)$.
We have $E^\perp=\gamma_2E$, and $\gamma_1$ preserves both $E$ and $E^\perp$, while $\gamma_2$ and $\gamma_3$ interchange them.
$d^\nabla\Psi=0$ and $\Pi(\Psi)=0$ imply $T_{ai}=T_{ia}$ and $\sum_i\gamma_iT_{ai}=0$.
The constraints give $T_{13}=-\gamma_2P+\gamma_1A,\ T_{23}=-\gamma_2A+\gamma_1B,\ T_{33}=-\gamma_2T_{13}+\gamma_1T_{23}=-P-B$.
Thus, $T$ is determined linearly by its independent components $(P,A,B)$.
Using the orthogonal decomposition $A=A_E+A_\perp$, we split these components as $(P,A,B)=(0,A_E,B)+(P,A_\perp,0)$.
Applying this decomposition for each component of $T$, we obtain $T=T_0+\Delta$, where $T_0$ and $\Delta$ are determined by $(0,A_E,B)$ and $(P,A_\perp,0)$, and given by 
\begin{align*}
T_0&=
\begin{pmatrix}
0&A_E&\gamma_1A_E\\
A_E&B&D_0\\
\gamma_1A_E&D_0&-B
\end{pmatrix},\quad
\Delta=
\begin{pmatrix}
P&A_\perp&-\gamma_2P+\gamma_1A_\perp\\
A_\perp&0&-\gamma_2A_\perp\\
-\gamma_2P+\gamma_1A_\perp&-\gamma_2A_\perp&-P
\end{pmatrix},\\
D_0&=-\gamma_2A_E+\gamma_1B.
\end{align*}
We have the estimate $|\Delta|^2=2|P|^2+4|A_\perp|^2+2|-\gamma_2P+\gamma_1A_\perp|^2\le6|P|^2+8|A_\perp|^2$ and $|T_0|\le |T|+|\Delta|\le C_1|T|$.
We now compute $H(\psi_0,T_0)$.
A direct calculation gives $H(\psi_0,T_0)=A_E\wedge B\wedge\gamma_1\Phi-\gamma_1A_E\wedge\Phi\wedge B+A_E\wedge\Phi\wedge D_0+\gamma_1A_E\wedge D_0\wedge\gamma_1\Phi$.
Writing $A_E=u\Phi+v\gamma_1\Phi$ with $u,v\in\mathbb R$, we obtain $A_E\wedge B\wedge\gamma_1\Phi-\gamma_1A_E\wedge\Phi\wedge B=-u\Phi\wedge\gamma_1\Phi\wedge B+u\Phi\wedge\gamma_1\Phi\wedge B=0$ and $A_E\wedge\Phi\wedge D_0+\gamma_1A_E\wedge D_0\wedge\gamma_1\Phi=-v\Phi\wedge\gamma_1\Phi\wedge D_0+v\Phi\wedge\gamma_1\Phi\wedge D_0=0$.
Therefore, we have $H(\psi_0,T_0)=0$.
Using $T=T_0+\Delta$, we obtain $H(\psi_0,T)-H(\psi_0,T_0)=\sum_a\Bigl[\Delta_{a1}\wedge T_{a2}\wedge\gamma_1\Phi+(T_0)_{a1}\wedge\Delta_{a2}\wedge\gamma_1\Phi+\Delta_{a1}\wedge\Phi\wedge T_{a3}+(T_0)_{a1}\wedge\Phi\wedge\Delta_{a3}\Bigr]$.
Therefore, we have $|H(\psi_0,T)|=|H(\psi_0,T)-H(\psi_0,T_0)|\le\sqrt2\,|\Phi|(|T|+|T_0|)|\Delta|\le C_2|\Phi||T|(|P|+|A_\perp|)$.
Combining this with~\eqref{eq:H-difference1}, we obtain \eqref{eq:H-component-bound}.
\end{proof}

\paragraph{Step 4}
The estimate~\eqref{eq:main} holds on $M^\circ$ with a  constant $C$.
\begin{proof}
Set $V_a=*_\Sigma(T_{a1}\wedge\Phi\wedge\gamma_1\Phi)$.
We have $|V_a|^2=|\Phi|^4|(T_{a1})_\perp|^2$, and summing over $a$ gives $|V|^2=|\Phi|^4\bigl(|P_\perp|^2+|A_\perp|^2+|(T_{31})_\perp|^2\bigr)$.
Since $T_{ai}=T_{ia}$ and $\sum_i\gamma_iT_{ai}=0$, we have $0=\operatorname{pr}_E\left(\sum_i\gamma_iT_{i1}\right)=\gamma_1P_E+\gamma_2A_\perp+\gamma_3(T_{31})_\perp$.
Thus, $|P_E|\le|A_\perp|+|(T_{31})_\perp|$.
It follows that
\begin{align}
|P|+|A_\perp|
&=\sqrt{|P_E|^2+|P_\perp|^2}+|A_\perp|
\le\sqrt{|P_\perp|^2+\bigl(|A_\perp|+|(T_{31})_\perp|\bigr)^2}+|A_\perp|\nonumber\\
&\le C_3\bigl(|P_\perp|^2+|A_\perp|^2+|(T_{31})_\perp|^2\bigr)^{1/2}
=C_3\,|\Phi|^{-2}|V|.\label{eq:PA-V}
\end{align}
To estimate $|P|+|A_\perp|$, we compare $V$ with $\nabla N$.
The Hodge star on the real spinor bundle and the volume form on $M$ are parallel, so we have $\nabla_{e_a}N=*_\Sigma\bigl(
T_{a1}\wedge\psi_2\wedge\psi_3+\psi_1\wedge T_{a2}\wedge\psi_3+\psi_1\wedge\psi_2\wedge T_{a3}\bigr)$.
Substituting~\eqref{eq:normalform} gives $\nabla_{e_a}N-V_a
= *_\Sigma\Bigl[tT_{a1}\wedge\Phi\wedge\gamma_1\Phi+(1+t)\psi_1\wedge T_{a2}\wedge\gamma_1\Phi+\psi_1\wedge\Phi\wedge T_{a3}\Bigr]$.
Using $t|\Phi|=|\psi_1|$ and $1+t\le2$, we obtain $|\nabla N-V|^2\le6|\Phi|^2|\psi_1|^2\sum_a\bigl(|T_{a1}|^2+|T_{a2}|^2+|T_{a3}|^2\bigr)=6|\Phi|^2|\psi_1|^2|T|^2$.
Using~\eqref{eq:PA-V} and the above bound on $|\nabla N-V|^2$, we obtain 
\begin{align}
|P|+|A_\perp|\le C_3|\Phi|^{-2}\bigl(|\nabla N|+|\nabla N-V|\bigr)\le C_3|\Phi|^{-2}|\nabla N|+C_4|\Phi|^{-1}|\psi_1||T|. \label{eq:PA-final}
\end{align}
Substituting~\eqref{eq:PA-final} into~\eqref{eq:H-component-bound} gives
\begin{align*}
|H(\psi,T)|
&\le C_5\left[
|\psi_1||T|^2
+|\Phi||T|\left(
\frac{|\nabla N|}{|\Phi|^2}
+\frac{|\psi_1||T|}{|\Phi|}
\right)\right]
\le C_6\left(\frac{|T|}{|\Phi|}|\nabla N|+|\psi_1||T|^2\right).
\end{align*}
Since $|\psi|=|\Psi|$ and $|T|=|\nabla\Psi|$, \eqref{eq:scale} gives $\frac{|T|}{|\Phi|}|\nabla N|\le\sqrt6\,\frac{|\nabla\Psi|}{|\Psi|}|\nabla N|$ and $|\psi_1||T|^2\le\frac{|N|}{|\Phi|^2}|T|^2\le6\,\frac{|\nabla\Psi|^2}{|\Psi|^2}|N|$.
Consequently, we obtain
\[
|H(\psi,T)|
\le C_7\left(
\frac{|\nabla\Psi|}{|\Psi|}|\nabla N|
+\frac{|\nabla\Psi|^2}{|\Psi|^2}|N|
\right).
\]
Substituting this into the equation of Step1 gives~\eqref{eq:main} at $q$.
Since $q\in M^\circ$ was arbitrary and $C$ is a universal constant, the estimate holds throughout $M^\circ$.
\end{proof}

\section{Proof of Lemma \ref{lem:geom}}\label{app:geom}
\paragraph{Step 1}
A geodesic generated by $\xi$ does not reach $Z$ at a finite arc length parameter.
\begin{proof}
Let $\gamma(t)$ be a unit-speed geodesic generated by $\xi$, and suppose that $\gamma(t)\in M^\circ$ for $0\le t<T<\infty$ and $\gamma(t)\to p\in Z$ as $t\uparrow T$.
Take $0\ne v\perp\xi$ at $\gamma(0)$.
Let $\Phi_t$ be the local flow of $\xi$, and let $J(t):=(\dd\Phi_t)_{\gamma(0)}v$ be its variation field.
Since $[\xi,J]=0$ and $\nabla_\xi\xi=0$, the field $J$ is a Jacobi field along $\gamma$, and $\frac{\dd}{\dd t}g(J,\xi)=g(\nabla_\xi J,\xi)=g(\nabla_J\xi,\xi)=0$, thus $J\perp\xi$ is preserved.
The structure equation \eqref{eq:xi-connection} and the assumption $\xi(\log f)=-\beta$ give $\frac{\dd}{\dd t}|J|^2=2g(\nabla_J\xi,J) =2\beta|J|^2-2\alpha g(\phi J,J)=2\beta|J|^2$.
Hence $\frac{\dd}{\dd t}\log\bigl(f(\gamma(t))|J(t)|\bigr)=-\beta+\frac{2\beta|J(t)|^2}{2|J(t)|^2}=0$.
Therefore $f(\gamma(t))|J(t)|=f(\gamma(0))|v|>0,\ (0\le t<T)$.
On the other hand, $\gamma$ is a unit-speed geodesic for the smooth metric $g$ and converges to $p\in M$ as $t\uparrow T$, so it extends smoothly to $t=T$.
The Jacobi equation is also a linear ordinary differential equation with smooth coefficients, so $J$ extends to $t=T$ and, in particular, is bounded.
Since $f$ is continuous on $M$ and $f(p)=0$, we have $f(\gamma(t))|J(t)|\to0$ as $t\uparrow T$, which contradicts $f(\gamma(T))|J(T)|=f(\gamma(0))|v|>0$.
\end{proof}

Since a three-dimensional orientable manifold is parallelizable, fix a smooth orthonormal frame $\{e_1,e_2,e_3\}$ on $M$.
For $x\in M^\circ$, take a local unit section $\xi$ of $L^\circ$ defined near $x$, and define
\[
 \ell=[\xi]:M^\circ\longrightarrow\mathbb{RP}^2,\qquad
 \ell(x):=[g_x(\xi_x,e_1(x)):g_x(\xi_x,e_2(x)):g_x(\xi_x,e_3(x))].
\]
Changing the local unit section only changes the sign of $\xi$ at each point, so $\ell(x)$ is independent of this choice.

\paragraph{Step 2}
The map $\ell$ extends continuously to $M$.
\begin{proof}
For \(p\in Z\), define $C_p:=\bigcap_{r>0}\overline{\ell(B_r(p)\setminus Z)}\subset\mathbb{RP}^2$.
Here $B_r(p)$ is the ball of radius $r$ centered at $p$.
The assumption $\mathcal H^1_{\loc}(Z)<\infty$ implies that $Z$ has no interior points, so $M^\circ$ is dense in $M$.
Thus $\overline{\ell(B_r(p)\setminus Z)}$ is nonempty and is compact by the compactness of $\mathbb{RP}^2$.
These sets decrease as $r\downarrow0$, so their intersection $C_p$ is also nonempty.

Next, take $r_0>0$ such that $B_{r_0}(p)$ is a normal ball.
By the compactness of the set of unit initial velocities at $p$ and the smoothness of the geodesic flow, we can fix a sufficiently small $0<2\varepsilon<r_0$ so that every unit-speed geodesic starting at $q\in B_\varepsilon(p)$ exists on $[-\varepsilon,\varepsilon]$.
For $\ell_0\in C_p$, take $q_j\in M^\circ$ such that $q_j\to p$ and $\ell(q_j)\to\ell_0$.
If $v\in T_pM$ is a unit vector representing $\ell_0$ at $p$, we can choose the signs of unit vectors $v_j$ in $L^\circ_{q_j}$ so that $v_j\to v$.
Here convergence of vectors means convergence of their components in the fixed frame.
Except for finitely many terms, $q_j\in B_\varepsilon(p)$, and the smooth dependence of the geodesic flow on initial values gives $\gamma_j\to\gamma$ in $C^1([-\varepsilon,\varepsilon])$ for $\gamma_j(t):=\exp_{q_j}(tv_j),\ \gamma(t):=\exp_p(tv)$.
Each $\gamma_j$ is a geodesic generated by $\xi$ in $M^\circ$.
If there is a point where the limit curve $\gamma$ enters $M^\circ$, the convergence of velocities and the continuity of $L^\circ$ imply that $\dot\gamma$ also belongs to $L^\circ$, and thus the geodesic reaches $Z$ since $\gamma(0)=p\in Z$, contradicting Step 1, so $\gamma([-\varepsilon,\varepsilon])\subset Z$.

Take distinct directions $\ell_1,\ldots,\ell_N\in C_p$, and let $v_1,\ldots,v_N$ be unit vectors representing them at $p$.
Since $\exp_p$ is injective in the normal ball, the arcs $\exp_p([-\varepsilon,\varepsilon]v_i)$ intersect each other only at $p$, and each arc has length $2\varepsilon$.
Thus $2\varepsilon N=\sum_{i=1}^N\mathcal H^1\bigl(\exp_p([-\varepsilon,\varepsilon]v_i)\bigr)=\mathcal H^1\left(\bigcup_{i=1}^N\exp_p([-\varepsilon,\varepsilon]v_i)\right)\le\mathcal H^1\bigl(Z\cap\overline{B_\varepsilon(p)}\bigr)<\infty$, so $C_p$ is a finite set.
Since \(\mathcal H^1_{\mathrm{loc}}(Z)<\infty\), the set $B_r(p)\setminus Z$ is connected for small $r$.
Thus $\overline{\ell(B_r(p)\setminus Z)}$ is a connected compact set, and the intersection \(C_p\) of this decreasing sequence is also connected.
Since $C_p$ is nonempty, finite, and connected, we have $C_p=\{\ell_p\}$.
Defining $\ell(p):=\ell_p$ extends $\ell$ to all of $M$.
\end{proof}

\paragraph{Step 3}
In a sufficiently small neighborhood of \(p\in Z\), the set \(Z\) is a single smooth geodesic arc.
In particular, we can take a suitable transverse disk \(\Sigma\) through \(p\) such that $\Sigma\cap Z=\{p\}$.
\begin{proof}
The geodesic through each \(p\in Z\) in the extended direction \(\ell_p\) lies in \(Z\).
Let $\gamma(t)$ be this geodesic and let $\gamma_j(t)$ be the approximating geodesics constructed in Step 2.
The convergence $\gamma_j\to\gamma$ is in $C^1$, and each $\gamma_j$ is tangent to $L^\circ$, so the continuity of $\ell$ gives $\lim_{j\to\infty}\ell(\gamma_j(t))=\ell(\gamma(t))$ at each time $t$.
Hence this geodesic is tangent to the extended line distribution at every point.
Thus \(Z\) is a union of extended geodesics.
Take a small smooth disk \(\Sigma\) through \(p\) transverse to \(\ell_p\).
By the continuity of the direction field, nearby geodesic leaves generated by $\xi$ are also uniformly transverse to \(\Sigma\).
If \(\Sigma\cap Z\) were not locally finite, each of its points would give a geodesic arc in the zero set, with these arcs mutually disjoint and of a common positive length, and we would have $\mathcal H^1(Z\cap K)=\infty$, contradicting the assumption.
Shrinking \(\Sigma\) further, we can arrange that $\Sigma\cap Z=\{p\}$.
Therefore, in a sufficiently small neighborhood of \(p\), the set \(Z\) is a single smooth geodesic arc.
\end{proof}

In what follows, choose a continuous local unit section $\xi$ of the line distribution extended in Step 2, and set $\phi X=\xi\times X$.
Take the transverse disk $\Sigma$ in Step 3 so that $T_p\Sigma=\xi(p)^\perp$, and choose its unit normal so that $n_\Sigma(p)=\xi(p)$.
After shrinking the disk suitably, we can arrange that $g(\xi,n_\Sigma)\ge\delta$ for some $\delta>0$.

\paragraph{Step 4}
Write $\xi=a e_1+b e_2+c e_3$ on the transverse disk $\Sigma$.
The functions $a,b$ satisfy the following equations on $\Sigma\setminus\{p\}$
\begin{align}
 \operatorname{div}_0(\mathsf A_a Da)&=\operatorname{div}_0G_a,\label{eq:div-a}\\
 \operatorname{div}_0(\mathsf A_b Db)&=\operatorname{div}_0G_b.\label{eq:div-b}
\end{align}
Here $D$ and $\operatorname{div}_0$ are the Euclidean gradient and divergence with respect to suitable coordinates on $\Sigma$.
The coefficients $\mathsf A_i$ are bounded on $\Sigma$, and $\mathsf A_i,\ G_i$ depend smoothly on $(x,y,a,b)$.
These equations are also uniformly elliptic on $\Sigma$.
\begin{proof}
Choose isothermal coordinates $(x,y)$ for $g|_\Sigma$ on $\Sigma$ so that $g|_\Sigma=e^{2\lambda}(dx^2+dy^2)$.
Take the frame $e_1=e^{-\lambda}\partial_x,\ e_2=e^{-\lambda}\partial_y,\ e_3=n_\Sigma$.
Write $\xi=a e_1+b e_2+c e_3,\ c=\sqrt{1-a^2-b^2}\ge\delta>0$.
We may assume that $a(p)=b(p)=0$.
By \eqref{eq:xi-connection}, we have $g(\nabla_x\xi,e_1)=e^\lambda\beta(1-a^2),\ g(\nabla_y\xi,e_1)+g(\nabla_x\xi,e_2)=-2e^\lambda\beta ab=m\,g(\nabla_x\xi,e_1)$ and $g(\nabla_y\xi,e_2)=e^\lambda\beta(1-b^2)=n\,g(\nabla_x\xi,e_1)$.
Setting $m=-\frac{2ab}{1-a^2},\ n=\frac{1-b^2}{1-a^2}$ and $F_1=g(\xi,\nabla_y e_1)-m\,g(\xi,\nabla_x e_1)+g(\xi,\nabla_x e_2),\ F_2=g(\xi,\nabla_y e_2)-n\,g(\xi,\nabla_x e_1)$, and using $g(\nabla_i\xi,e_j)=\partial_i g(\xi,e_j)-g(\xi,\nabla_i e_j)$, we obtain
\begin{equation}\label{eq:real-first-order}
 \begin{aligned}
 a_y&=ma_x-b_x+F_1,\\
 b_y&=na_x+F_2.
 \end{aligned}
\end{equation}
The functions $F_1,F_2$ are smooth functions of $(x,y,a,b)$ and contain no derivatives of $a,b$.
Let $D$ and $\operatorname{div}_0$ be the Euclidean gradient and divergence with respect to the coordinates $(x,y)$, and set
\begin{equation}\label{eq:AaGa}
 \mathsf A_a=\begin{pmatrix}n&0\\-m&1\end{pmatrix},\qquad
 G_a=\binom{-F_2}{F_1},
\end{equation}
\begin{equation}\label{eq:AbGb}
 \mathsf A_b=\begin{pmatrix}1&-m/n\\0&1/n\end{pmatrix},\qquad
 G_b=\binom{F_1-mF_2/n}{F_2/n}.
\end{equation}
Then \eqref{eq:real-first-order} gives
\begin{align}
 \mathsf A_a Da-G_a
 &=\binom{na_x+F_2}{-ma_x+a_y-F_1}
 =\binom{b_y}{-b_x},\label{eq:flux-a}\\
 \mathsf A_b Db-G_b
 &=\binom{b_x-(m/n)b_y-F_1+mF_2/n}{(b_y-F_2)/n}
 =\binom{-a_y}{a_x}.\label{eq:flux-b}
\end{align}
Since the divergence of the right-hand sides is zero, we obtain $\operatorname{div}_0(\mathsf A_a Da)=\operatorname{div}_0G_a$ and $\operatorname{div}_0(\mathsf A_b Db)=\operatorname{div}_0G_b$.

We compute the principal symbols.
We have $\operatorname{div}_0(\mathsf A_a Dv)=\partial_x(nv_x)+\partial_y(-mv_x+v_y)=nv_{xx}-mv_{xy}+v_{yy}+(n_x-m_y)v_x$ and $\operatorname{div}_0(\mathsf A_b Dv)=\partial_x\left(v_x-\frac mn v_y\right)+\partial_y\left(\frac1n v_y\right)=v_{xx}-\frac mn v_{xy}+\frac1n v_{yy}   +\left\{\partial_y\left(\frac1n\right)           -\partial_x\left(\frac mn\right)\right\}v_y$.
Since $\operatorname{div}_0G_a,\operatorname{div}_0G_b$ also contain only derivatives of $a,b$ of order at most one, the principal symbols are
\begin{align*}
 \sigma_a(\tau)
 &=n\tau_x^2-m\tau_x\tau_y+\tau_y^2
 =\frac{c^2|\tau|^2+(a\tau_x+b\tau_y)^2}{1-a^2}=\tau^T\mathsf A_a\tau,\\
 \sigma_b(\tau)
 &=\tau_x^2-\frac mn\tau_x\tau_y+\frac1n\tau_y^2
 =\frac{c^2|\tau|^2+(a\tau_x+b\tau_y)^2}{1-b^2}=\tau^T\mathsf A_b\tau.
\end{align*}
In general, uniform ellipticity of a differential operator \(L\) of order \(m\) means that its principal symbol is invertible and there exists a constant \(C\), independent of the point \(x\) in the domain, such that $\|\sigma_L(x,\tau)^{-1}\|\le C|\tau|^{-m},\ (\tau\ne0)$.
We show that these principal symbols are uniformly elliptic.
Noting that $c^2|\tau|^2\le c^2|\tau|^2+(a\tau_x+b\tau_y)^2\le(c^2+a^2+b^2)|\tau|^2=|\tau|^2$ and $\delta^2\le c^2\le1-a^2,1-b^2\le1$, we obtain $\delta^2|\tau|^2\le\sigma_j(\tau)=\tau^T\mathsf A_j\tau \le\delta^{-2}|\tau|^2,\ (j=a,b)$.
\end{proof}

In the setting of Step 4, we prove the following lemma to show that $a,b$ extend smoothly to the point $p$.
\begin{lemma}\label{lem:real-removal}
Let $\Omega\subset\R^2$ be an open disk containing the origin, and set $\Omega^\ast:=\Omega\setminus\{0\}$.
Let $U\subset\R^2$ be an open set, and let $\mathsf A_a,\mathsf A_b\in C^\infty(\Omega\times U,\R^{2\times2}),\ G_a,G_b\in C^\infty(\Omega\times U,\R^2)$.
Suppose that the real-valued functions $a,b\in C^0(\Omega)\cap C^\infty(\Omega^\ast)$ satisfy $(a(x,y),b(x,y))\in U$ at every $(x,y)\in\Omega$ and that, for some constants $0<\lambda_0\le\Lambda_0$, $ \tau^T\mathsf A_j\tau\ge\lambda_0|\tau|^2,\ \|\mathsf A_j\|\le\Lambda_0,\ \tau\in\R^2,\ j=a,b$ holds throughout $\Omega$.
Let $D$ and $\operatorname{div}_0$ denote the Euclidean gradient and divergence on $\Omega$, and assume that \eqref{eq:div-a} and~\eqref{eq:div-b} hold on $\Omega^\ast$.
Then $a,b\in C^\infty(\Omega)$, and these equations also hold on $\Omega$.
\end{lemma}
\begin{proof}
The smoothness of the coefficients and the continuity of $(a,b)$ imply that $G_a,G_b$ are locally bounded.
All integrals are with respect to $\dd x\dd y$.
First, we show that $a,b\in W^{1,2}_{\loc}(\Omega)$ and that \eqref{eq:div-a} and~\eqref{eq:div-b} hold weakly throughout $\Omega$.
Write $v=a$ or $b$, and denote the corresponding coefficient and right-hand side by $\mathsf A,G$.
On $\Omega^\ast$, we have $\operatorname{div}_0(\mathsf A Dv)=\operatorname{div}_0G$.
For $\chi\in C_c^\infty(\Omega^\ast,\R)$, using $\chi^2v$ as a test function gives $\int\mathsf A Dv\cdot D(\chi^2v)=-\int\operatorname{div}_0(\mathsf A Dv)\,\chi^2v=-\int\operatorname{div}_0G\,\chi^2v=\int G\cdot D(\chi^2v)$, and hence
\begin{align}
 \int_{\Omega^\ast}\chi^2\mathsf A Dv\cdot Dv
 &=\int_{\Omega^\ast}\chi^2G\cdot Dv
   +2\int_{\Omega^\ast}\chi v\,G\cdot D\chi
   -2\int_{\Omega^\ast}\chi v\,\mathsf A Dv\cdot D\chi.
 \label{eq:energy-test}
\end{align}
By assumption, $\mathsf A Dv\cdot Dv\ge\lambda_0|Dv|^2$ and $|\mathsf A Dv|\le\Lambda_0|Dv|$, so \eqref{eq:energy-test} gives
\begin{align*}
 \lambda_0\int_{\Omega^\ast}\chi^2|Dv|^2
 &\le\int_{\Omega^\ast}\chi^2|G||Dv|
   +2\int_{\Omega^\ast}|\chi v|\,|G||D\chi|
   +2\Lambda_0\int_{\Omega^\ast}|\chi v|\,|Dv||D\chi|.
\end{align*}
Applying Young's inequality to each term on the right-hand side gives $\chi^2|G||Dv|\le\frac{\lambda_0}{4}\chi^2|Dv|^2+\frac1{\lambda_0}\chi^2|G|^2$, $2|\chi v|\,|G||D\chi|\le\chi^2|G|^2+v^2|D\chi|^2$, and $2\Lambda_0|\chi v|\,|Dv||D\chi|\le\frac{\lambda_0}{4}\chi^2|Dv|^2+\frac{4\Lambda_0^2}{\lambda_0}v^2|D\chi|^2$, and hence
\begin{gather}
 \lambda_0\int_{\Omega^\ast}\chi^2|Dv|^2
 \!\le\!\frac{\lambda_0}{2}\int_{\Omega^\ast}\chi^2|Dv|^2
   +\left(1+\frac1{\lambda_0}\right)\int_{\Omega^\ast}\chi^2|G|^2
   +\left(1+\frac{4\Lambda_0^2}{\lambda_0}\right)\int_{\Omega^\ast} v^2|D\chi|^2.\nonumber\\
 \frac{\lambda_0}{2}\int_{\Omega^\ast}\chi^2|Dv|^2
 \le\left(1+\frac1{\lambda_0}\right)\int_{\Omega^\ast}\chi^2|G|^2
   +\left(1+\frac{4\Lambda_0^2}{\lambda_0}\right)\int_{\Omega^\ast} v^2|D\chi|^2.\nonumber\\
  \int_{\Omega^\ast}\chi^2|Dv|^2
 \le C\int_{\Omega^\ast}\bigl(v^2|D\chi|^2+\chi^2|G|^2\bigr)\label{eq:energy-pair}
\end{gather}
Here $C$ depends only on $\lambda_0,\Lambda_0$.
Fix $\theta\in C_c^\infty(\Omega)$ that equals $1$ near the origin.
Choose a smooth function $\rho_\varepsilon\in C^\infty(\R^2)$ that equals $0$ near the origin and $1$ on $\R^2\setminus B_\varepsilon(0)$ such that $0\le\rho_\varepsilon\le1$ and $\int_{\R^2}|D\rho_\varepsilon|^2\to0,\ (\varepsilon\to0)$.
Take the cutoff function $\chi=\theta\rho_\varepsilon\in C_c^\infty(\Omega^\ast)$.
By taking $\varepsilon$ sufficiently small, we may assume that $\|D\rho_\varepsilon\|_{L^2(\R^2)}^2\le1$.
First, we have $\int_\Omega|D\chi|^2=\int_\Omega|\rho_\varepsilon D\theta+\theta D\rho_\varepsilon|^2 \le2\int_\Omega|D\theta|^2+2\|\theta\|_\infty^2$ and $ \int_\Omega\chi^2=\int_\Omega\theta^2\rho_\varepsilon^2 \le\int_\Omega\theta^2$.
Thus, for $v=a,b$, \eqref{eq:energy-pair} gives
\begin{align*}
 \int_\Omega|D(\chi v)|^2
 &=\int_{\Omega^\ast}|\chi Dv+vD\chi|^2
 \le2\int_{\Omega^\ast}\chi^2|Dv|^2
   +2\int_{\Omega^\ast}v^2|D\chi|^2\\
 &\le2\int_{\Omega^\ast}\chi^2|Dv|^2
   +2\|v\|_{L^\infty(\operatorname{supp}\theta)}^2\int_\Omega|D\chi|^2\\
 &\le2C\int_{\Omega^\ast}\bigl(v^2|D\chi|^2+\chi^2|G|^2\bigr)
   +2\|v\|_{L^\infty(\operatorname{supp}\theta)}^2\int_\Omega|D\chi|^2\\
 &\le2(C+1)\|v\|_{L^\infty(\operatorname{supp}\theta)}^2\int_\Omega|D\chi|^2
   +2C\|G\|_{L^\infty(\operatorname{supp}\theta)}^2\int_\Omega\chi^2
 \le C_\theta.
\end{align*}
Here $C_\theta:=4(C+1)\|v\|_{L^\infty(\operatorname{supp}\theta)}^2\left(\int_\Omega|D\theta|^2+\|\theta\|_\infty^2\right)+2C\|G\|_{L^\infty(\operatorname{supp}\theta)}^2\int_\Omega\theta^2$.
The function $\chi v$ belongs to $W^{1,2}$ on all of $\Omega$, and $\chi v\to\theta v,\ (\varepsilon\to0)$ in $L^2$.
Weak compactness gives $\theta v\in W^{1,2}(\Omega)$.
Since they are already smooth away from the origin, we obtain $a,b\in W^{1,2}_{\loc}(\Omega)$.

Next, we extend the weak equation to the origin.
Take $\varphi\in C_c^\infty(\Omega)$.
Using $\varphi\rho_\varepsilon$ as a test function in the equation on $\Omega^*$ gives $0=\int_{\Omega^\ast}(\mathsf A Dv-G)\cdot D(\varphi\rho_\varepsilon)=\int_{\Omega^\ast}\rho_\varepsilon(\mathsf A Dv-G)\cdot D\varphi +\int_{\Omega^\ast}\varphi(\mathsf A Dv-G)\cdot D\rho_\varepsilon$.
Since $\mathsf A Dv-G\in L^2(\operatorname{supp}\varphi)$, we have $\left|\int_{\Omega^\ast}\varphi(\mathsf A Dv-G)\cdot D\rho_\varepsilon\right|\le\|\varphi\|_\infty\|\mathsf A Dv-G\|_{L^2(\operatorname{supp}\varphi)} \|D\rho_\varepsilon\|_{L^2(\R^2)}\longrightarrow0$.
For the remaining term $\int_{\Omega^\ast}\rho_\varepsilon(\mathsf A Dv-G)\cdot D\varphi$, we can pass to the limit by the dominated convergence theorem, obtaining
\[
 \int_\Omega\mathsf A Dv\cdot D\varphi=\int_\Omega G\cdot D\varphi.
\]
Thus the equation extends to all of $\Omega$.

Finally, we show that $a,b\in C^\infty(\Omega)$ and that \eqref{eq:div-a} and~\eqref{eq:div-b} hold throughout $\Omega$.
By the continuity of $(a,b)$, $\mathsf A_a,\mathsf A_b$ are continuous, uniformly elliptic, and bounded, including at the origin.
Also, $G_a,G_b$ are locally bounded.
Applying the interior H\"older estimate for weak solutions of scalar divergence-form equations with continuous coefficients~\cite[Theorem~6.1]{Vita2025} to each component gives $a,b\in C^{0,\gamma}_{\loc}(\Omega)$ for some $\gamma\in(0,1)$.
Since the coefficients and right-hand sides depend smoothly on $(x,y,a,b)$, we have $\mathsf A_j,G_j\in C^{0,\gamma}_{\loc}(\Omega)$.
Thus the divergence-form Schauder estimate~\cite[Theorem~6.2]{Vita2025} gives $a,b\in C^{1,\gamma}_{\loc}(\Omega)$.
Furthermore, by repeated application of~\cite[Theorem~6.3]{Vita2025}, we obtain $(a,b)\in C^{k,\gamma}_{\loc}(\Omega)\Rightarrow \mathsf A_j,G_j\in C^{k,\gamma}_{\loc}(\Omega)\Rightarrow (a,b)\in C^{k+1,\gamma}_{\loc}(\Omega),\ k\ge1$.
Therefore, $a,b\in C^\infty(\Omega)$.
\end{proof}

\paragraph{Step 5}
$L^\circ$ extends uniquely to a smooth real line bundle $L$ on all of $M$, and each of its local unit sections defines a trans-Sasakian $\eta$-Einstein structure.
\begin{proof}
The functions $(a,b)$ in Step 4 satisfy the assumptions of Lemma~\ref{lem:real-removal}.
The restriction $\xi|_\Sigma$ is smooth, including at $p$.
With $\xi|_\Sigma$ as the initial velocity, set $F(t,q)=\exp_q\bigl(t\xi(q)\bigr),\ q\in\Sigma$.
The differential $dF_{(0,p)}(s,v)=s\xi(p)+v$ is an isomorphism, so $F$ is a local diffeomorphism near $(0,p)$.
For $q\ne p$, $F(t,q)$ coincides with a geodesic leaf of the original $\xi$, so $F_*\partial_t$ gives a smooth extension of $\xi$.
Thus $L=\R\xi$ is smooth near $p$.
For a smooth local orthonormal frame $\{e_i\}$, the type of the structure extends smoothly by $ \beta=\frac12\diver\xi,\ \alpha=-\frac12\sum_{i=1}^3 g(\nabla_{e_i}\xi,\phi e_i)$, and \eqref{eq:xi-connection} continues to hold after extension by the density of $M^\circ$.
Similarly, taking the trace and the $(\xi,\xi)$ component of the $\eta$-Einstein condition $\Ric=\lambda g+\nu\eta\otimes\eta$ on $M^\circ$ gives $\Scal=3\lambda+\nu$ and $\Ric(\xi,\xi)=\lambda+\nu$.
The formulas $\lambda=\frac12\{\Scal-\Ric(\xi,\xi)\}$ and $\nu=\frac12\{3\Ric(\xi,\xi)-\Scal\}$ therefore define smooth extensions of $\lambda$ and $\nu$.
The smooth tensor $\Ric-\lambda g-\nu\eta\otimes\eta$ vanishes on the dense set $M^\circ$, so it vanishes throughout the extension domain by continuity.
These local extensions agree on $M^\circ$, and the extension is unique by the density of $M^\circ$.
\end{proof}

\section{Transversal holomorphic geometry}\label{app:transversal holomorphic geometry}
We summarize the terminology from transversal holomorphic geometry used in this paper.
For a one-dimensional foliation $\mathcal F$ on a three-dimensional manifold, choose flow box coordinates $(t_i,z_i):U_i\longrightarrow I_i\times V_i$ near each point, where $I_i\subset\R$ is an open interval and $V_i\subset\C$ is an open set.
A curve with $z_i$ fixed is called a plaque.
If such coordinates can be chosen so that the transverse coordinate changes on overlaps can locally be written as $z_j=f_{ij}(z_i)$ with $f_{ij}$ locally biholomorphic, then $\mathcal F$ is called a \textcolor{blue}{\textbf{transversely holomorphic foliation}}~\cite[Section~2.1]{BiswasKasuya2023}.

Next, consider a smooth complex line bundle $E\to M$ over a three-dimensional manifold $M$ with a transversely holomorphic foliation $\mathcal F$.
Choose nowhere-vanishing smooth local sections $e_i$ of $E$ on the above flow boxes $U_i$ and suppose that, locally on overlaps, we can write $e_j=g_{ij}(z_i)e_i$ with $g_{ij}$ a nowhere-vanishing holomorphic function.
The structure defined by this system of trivializations is called a basic holomorphic structure on $E$, and $E$ equipped with this structure is called a \textcolor{blue}{\textbf{basic holomorphic line bundle}}~\cite[Definition~2.1]{BiswasKasuya2023}.
Writing a smooth section $s\in\Gamma(M,E)$ of $E$ as $s|_{U_i}=f_i e_i$, we call $s$ a \textcolor{blue}{\textbf{basic section}} if each $f_i$ is independent of $t_i$~\cite[Definition~2.2]{BiswasKasuya2023}.

Consider a flat partial connection $\mathcal D$ on $E$ along the complex distribution spanned by $\partial_{t_i}$ and $\partial_{\bar z_i}$ in the above flow box coordinates.
The statement that $\mathcal D$ defines a basic holomorphic line bundle structure on $E$ means that a basic holomorphic structure is obtained by choosing nowhere-vanishing $\mathcal D$-parallel frames $e_i\in\Gamma(U_i,E)$ on sufficiently small flow boxes $U_i$~\cite[Section~2.2]{BiswasKasuya2023}.
For the basic holomorphic structure defined by $\mathcal D$, using its $\mathcal D$-parallel frames $e_i$ gives $\mathcal D_{\partial_{t_i}}s=(\partial_{t_i}f_i)e_i+f_i\mathcal D_{\partial_{t_i}}e_i=(\partial_{t_i}f_i)e_i$, so $s$ is basic if and only if it is $\mathcal D$-parallel along the leaves.

\section{Proof of Proposition \ref{prop:distinct-ricci-eigenvalues}}\label{app:r3}
If $\Ric=fg$ holds, $f$ is constant by Schur's lemma.
In three dimension, then $g$ has constant curvature.
Thus there is a nonempty open set on which the Ricci endomorphism is not a scalar multiple.
Since $\Psi$ has rank three almost everywhere, this open set contains a point $q$ with $\rank_\R\Psi_q=3$.

Suppose that $\Ric^\sharp_q$ has a repeated eigenvalue $\lambda$, then its eigenvalues are $\lambda,\lambda,\mu$ with $\mu\ne\lambda$.
Choose a corresponding orthonormal eigenbasis $e_1,e_2,e_3$ of $T_qM$.
\eqref{eq:ricci-constraint} give $0=\cR(\Psi)_q-\lambda\Pi(\Psi)_q=\sum_{i=1}^3\gamma\bigl((\Ric^\sharp_q-\lambda\Id)e_i\bigr)\Psi_q(e_i)=(\mu-\lambda)\gamma(e_3)\Psi_q(e_3)$.
Since $\mu-\lambda\ne0$, we obtain $\Psi_q(e_3)=0$, contradicting $\rank_\R\Psi_q=3$.
Therefore, $\Ric^\sharp_q$ has three pairwise distinct eigenvalues.

Choose a coordinate neighborhood $U$ of $q$ on which $\Psi$ is nowhere zero and the Ricci eigenvalues $\rho_1,\rho_2,\rho_3$ remain pairwise distinct.
Choose a smooth orthonormal Ricci eigenframe $e_1,e_2,e_3$ and set $(c_1,c_2,c_3)=(\rho_2-\rho_3,\rho_3-\rho_1,\rho_1-\rho_2)$.
By $\Pi(\Psi)=0$ and~\eqref{eq:ricci-constraint}, we have $ \sum_{i=1}^3\gamma_i\Psi(e_i)=0$ and $\sum_{i=1}^3\rho_i\gamma_i\Psi(e_i)=0$.
Subtracting $\rho_3$ and $\rho_2$ times the first equation from the second, respectively, we obtain $0=\sum_{i=1}^3(\rho_i-\rho_3)\gamma_i\Psi(e_i)=-c_2\gamma_1\Psi(e_1)+c_1\gamma_2\Psi(e_2)$ and $0=\sum_{i=1}^3(\rho_i-\rho_2)\gamma_i\Psi(e_i)=c_3\gamma_1\Psi(e_1)-c_1\gamma_3\Psi(e_3)$.
Define the smooth spinor $\chi:=-c_1^{-1}\gamma_1\Psi(e_1)$ on $U$.
The preceding two equalities imply $\Psi(e_i)=-\frac{c_i}{c_1}\gamma_i\gamma_1\Psi(e_1)=c_i\gamma_i\chi$ for $i=1,2,3$.
Since $\Psi$ is nowhere zero on $U$, so is $\chi$.
Since $\Sigma M$ has complex rank two, $\chi,\mathfrak j\chi$ therefore form a smooth complex frame on $U$.
For any $\Phi\in\operatorname{RS}(M,g)$, the same argument gives a smooth spinor $\chi_\Phi$ on $U$ such that $\Phi(e_i)=c_i\gamma_i\chi_\Phi$ for $i=1,2,3$.
There are therefore unique functions $a,b\in C^\infty(U,\C)$ such that $\chi_\Phi=a\chi+b\mathfrak j\chi$.
Thus we have $\Phi(e_i)=c_i\gamma_i(a\chi+b\mathfrak j\chi) =a\Psi(e_i)+b(\mathfrak j\Psi)(e_i)$, and hence $\Phi=a\Psi+b\mathfrak j\Psi$ on $U$.

To show that $a$ and $b$ are constant, set $v_i=(e_i a)\chi+(e_i b)\mathfrak j\chi$ and $S=\sum_{i=1}^3\gamma_iv_i$.
We have $\Pi(\mathfrak j\Psi)=\mathfrak j(\Pi\Psi)=0$ and $D_{TM}(\mathfrak j\Psi)=\mathfrak j(D_{TM}\Psi)=0$, so $\mathfrak j\Psi$ is also an RS field.
Then we have $0=(D_{TM}\Phi)(e_j)=a(D_{TM}\Psi)(e_j)+b\bigl(D_{TM}(\mathfrak j\Psi)\bigr)(e_j)+\sum_{i=1}^3\gamma_i\bigl((e_i a)\Psi(e_j)+(e_i b)(\mathfrak j\Psi)(e_j)\bigr)=c_j\sum_{i=1}^3\gamma_i\gamma_j\bigl((e_i a)\chi+(e_i b)\mathfrak j\chi\bigr)=c_j\sum_{i=1}^3\gamma_i\gamma_jv_i =-\gamma_j\sum_{i=1}^3\gamma_iv_i-2v_j=-\gamma_jS-2v_j$.
It follows that $v_j=-\frac12\gamma_jS$, and then $S=\sum_{j=1}^3\gamma_jv_j=-\frac12\sum_{j=1}^3\gamma_j^2S=\frac32S$ holds.
Thus we have $S=0$ and $v_j=0$.
Therefore we obtain $\dd a=\dd b=0$.
Since $U$ is connected, $a$ and $b$ are constant on $U$.
The RS field $\Phi-a\Psi-b\mathfrak j\Psi$ vanishes on $U$, so unique continuation for the twisted Dirac operator~\cite{Bar1999} and the connectedness of $M$ imply that it vanishes on all of $M$.


\begin{thebibliography}{99}

\bibitem{RaritaSchwinger1941}
W.~Rarita and J.~Schwinger,
\emph{On a theory of particles with half-integral spin},
Phys. Rev. \textbf{60} (1941), 61.
\url{https://doi.org/10.1103/PhysRev.60.61}.

\bibitem{LawsonMichelsohn1989}
H.~B.~Lawson, Jr. and M.-L.~Michelsohn,
\emph{Spin Geometry},
Princeton Mathematical Series, vol.~38,
Princeton University Press, Princeton, NJ, 1989.
\url{https://doi.org/10.1515/9781400883912}.

\bibitem{HommaSemmelmann2019}
Y.~Homma and U.~Semmelmann,
\emph{The kernel of the Rarita-Schwinger operator on Riemannian spin manifolds},
Comm. Math. Phys. \textbf{370} (2019), no.~3, 853--871.
\url{https://doi.org/10.1007/s00220-019-03324-8}.

\bibitem{BarMazzeo2021}
C.~B\"ar and R.~Mazzeo,
\emph{Manifolds with many Rarita-Schwinger fields},
Comm. Math. Phys. \textbf{384} (2021), no.~1, 533--548.
\url{https://doi.org/10.1007/s00220-021-04030-0}.

\bibitem{OhnoTomihisa2023}
S.~Ohno and T.~Tomihisa,
\emph{Rarita-Schwinger fields on nearly K\"ahler manifolds},
Differential Geom. Appl. \textbf{91} (2023), Paper No.~102068.
\url{https://doi.org/10.1016/j.difgeo.2023.102068}.

\bibitem{Ohno2023G2}
S.~Ohno,
\emph{Rarita-Schwinger fields on nearly parallel $G_2$-manifolds},
J. Geom. Phys. \textbf{194} (2023), Paper No.~105024.
\url{https://doi.org/10.1016/j.geomphys.2023.105024}.

\bibitem{Nguyen2023}
M.~L.~Nguyen,
\emph{Pin(2)-equivariance property of the Rarita-Schwinger--Seiberg--Witten equations},
J. Geom. Anal. \textbf{33} (2023), no.~10, Paper No.~336.
\url{https://doi.org/10.1007/s12220-023-01404-9}.

\bibitem{Ohno2026LieGroups}
S.~Ohno,
\emph{Left-invariant Rarita-Schwinger fields on 3-dimensional Lie groups},
arXiv:2608.17392 (2026).
\url{https://arxiv.org/abs/2608.17392}.

\bibitem{Matsuno2025}
S.~Matsuno,
Three-dimensional almost contact metric manifolds revisited via the
Newman--Penrose formalism,
\emph{J. Geom.} \textcolor{blue}{\textbf{117}} (2026), Article~34.
\href{https://doi.org/10.1007/s00022-026-00816-6}{doi:10.1007/s00022-026-00816-6}.

\bibitem{Matsuno2026}
S.~Matsuno,
\emph{Three-dimensional noncompact trans-Sasakian $\eta$-Einstein manifolds},
arXiv:2610.39971 [math.DG] (2026).
\url{https://arxiv.org/abs/2610.39971}.

\bibitem{Mukud}
Prachi, M.-A.~Lawn and M.~M.~Tripathi,
\emph{Newman--Penrose formalism in $3$-dimensional trans-Sasakian manifolds},
arXiv:2605.20000 [math.DG] (2026).
\url{https://arxiv.org/abs/2605.20000}.

\bibitem{BairdWood1991}
P.~Baird and J.~C.~Wood,
\emph{Harmonic morphisms and conformal foliations by geodesics of three-dimensional space forms},
J. Austral. Math. Soc. (Series A) \textbf{51} (1991), 118--153.
\url{https://doi.org/10.1017/S1446788700033358}.

\bibitem{BiswasSchumacher}
I.~Biswas and G.~Schumacher,
\emph{Vector bundles on Sasakian manifolds},
Adv. Theor. Math. Phys. \textbf{14} (2010), no.~2, 541--561.
\url{https://doi.org/10.4310/ATMP.2010.v14.n2.a5}.
\url{https://arxiv.org/abs/0809.3892}.

\bibitem{Hawkins2008}
E.~Hawkins,
\emph{A groupoid approach to quantization},
J. Symplectic Geom. \textbf{6} (2008), no.~1, 61--125.
\url{https://doi.org/10.4310/JSG.2008.v6.n1.a4}.
\url{https://arxiv.org/abs/math/0612363}.

\bibitem{McMullenRS}
C.~McMullen,
\emph{Riemann Surfaces},
Course Notes, Harvard University, Math~213b.
\url{https://people.math.harvard.edu/~ctm/papers/home/text/class/harvard/213b/course/course.pdf}.

\bibitem{Bar1999}
C.~B\"ar,
\emph{Zero sets of solutions to semilinear elliptic systems of first order},
Invent. Math. \textbf{138} (1999), 183--202.
\href{https://doi.org/10.1007/s002220050346}{doi:10.1007/s002220050346}.
\href{https://arxiv.org/abs/math/9805065}{arXiv:math/9805065}.

\bibitem{Gerner2021}
W.~Gerner,
\emph{Typical field lines of Beltrami flows and boundary field line behaviour of Beltrami flows on simply connected, compact, smooth manifolds with boundary},
Ann. Global Anal. Geom. \textbf{60} (2021), 65--82.
\url{https://doi.org/10.1007/s10455-021-09768-3}.


\bibitem{BrandingEtAl}
V.~Branding, S.~Montaldo, C.~Oniciuc, A.~Ratto,
Unique continuation properties for polyharmonic maps between Riemannian manifolds,
\emph{Canad. J. Math.} \textbf{75} (2023), 1--28.





\bibitem{Vita2025}
S.~Vita,
\emph{Notes on Schauder estimates by scaling for second order linear elliptic PDEs in divergence form},
arXiv:2510.01765v1 (2025).
\href{https://arxiv.org/abs/2510.01765v1}{arXiv:2510.01765v1}.


\bibitem{BiswasKasuya2023}
I.~Biswas and H.~Kasuya,
\emph{Higgs bundles and flat connections over compact Sasakian manifolds, II: quasi-regular bundles},
arXiv:2110.10644v2 (2023).
\url{https://arxiv.org/abs/2110.10644v2}.

\bibitem{WangZhang2013}
G.~Wang and Y.~Zhang,
\emph{The Sasaki--Ricci flow on Sasakian 3-spheres},
Commun. Math. Stat. \textbf{1} (2013), no.~1, 43--71.
\url{https://doi.org/10.1007/s40304-013-0005-7}.
\url{https://arxiv.org/abs/1303.2491}.

\bibitem{Friedrich-Kim}
E.~C.~Kim and T.~Friedrich,
\emph{The Einstein--Dirac equation on Riemannian spin manifolds},
J. Geom. Phys. \textbf{33} (2000), no.~1--2, 128--172.
\url{https://doi.org/10.1016/S0393-0440(99)00043-1}.
\url{https://arxiv.org/abs/math/9905095}.

\end{thebibliography}
\end{document}